\documentclass[a4paper,10pt]{article} 
\usepackage{amsmath,amssymb,amsthm}
\usepackage[english]{babel}

\newtheorem{theo}{Theorem}[section] 
\newtheorem{defi}[theo]{Definition}
\newtheorem{lemm}[theo]{Lemma} 
\newtheorem{prop}[theo]{Proposition}
\newtheorem{coro}[theo]{Corollary}

\newcommand{\Ra}{\mathbb R}                   

\newcommand{\Ca}{\mathbb C}                   

\newcommand{\scal}[1]{\langle #1 \rangle}

\newcommand{\finpreuve}{\hfill $\Box$}

\newcommand{\name}{$\underline{\qquad \qquad}$} \newcommand{\etan}{$ \& $ }

\begin{document}

\author{  Jean-Marc Bouclet}
\title{Low frequency estimates and local energy decay for asymptotically euclidean Laplacians}

\maketitle

\begin{abstract} 
For riemannian metrics $ G $ on $ \Ra^d $ which are  long range perturbations of the flat one, we prove estimates for $(-\Delta_G-\lambda-i\epsilon)^{-n}$
as $ \lambda \rightarrow 0 $, which are uniform with respect to $ \epsilon $, for all $n \leq [d/2] + 1$ in odd dimension and $ n \leq d/2 $ in even dimension.
We also give applications to the time decay of Schr\"odinger and Wave (or Klein-Gordon) equations.  
\end{abstract}

\setcounter{section}{0}

\section{Introduction and results} \label{introduction}
\setcounter{equation}{0}
Let $ G = (G^{jk}) $ be a riemannian metric on $ \Ra^d $ which is asymptotically euclidean in the sense that, for some $ \rho > 0 $,
\begin{eqnarray}
|\partial^{\alpha} (G^{jk}(x)- \delta_{jk}) | \leq C_{\alpha} \scal{x}^{- \rho - |\alpha|} , \label{cometriquelongueportee}
\end{eqnarray}
$ \delta_{jk}$ being the Kronecker symbol. In other words, (the coefficients of) $ G - I $ belongs to the symbol class $ S^{-\rho} $ of functions such that $ |\partial^{\alpha}a(x)| \lesssim \scal{x}^{-\rho-|\alpha|} $. In the sequel we shall also refer to $ G$ as a long range metric. The Laplacian $ \Delta_G $ reads 
\begin{eqnarray}
 \Delta_G = \mbox{det}\ G(x)^{-1/2} \frac{\partial}{\partial x_j} \left( \mbox{det}\ G(x)^{1/2}  G_{jk}(x)  \frac{\partial}{\partial x_k} \right) , \label{LaplacienG}
\end{eqnarray} 
 using the summation convention as well as the standard notation $ (G_{jk}) := (G^{jk})^{-1} $, and is (formally) self-adjoint with respect to the measure
 $$ d_G x = \mbox{det} \ G(x)^{1/2} dx . $$
Since $ \mbox{det} \ G(x)^{1/2} $ is bounded from above and below, the spaces $ L^2 (\Ra^d, dx) $ and $ L^2 (\Ra^d, d_G x) $ coincide and have equivalent norms. We will thus use the unambiguous notation $  L^2 (\Ra^d) $ (or  $ L^2 $) in the sequel. By $ \Delta_G $ we will also denote the self-adjoint realization of (\ref{LaplacienG}), whose domain is $ H^2 $.

The main purpose of this paper is to investigate the low frequency estimates for powers of the resolvent of $ - \Delta_G $, namely the behaviour of
$$ (- \Delta_G - z)^{-n} \ \ \mbox{as} \ z \ \mbox{approaches} \ 0, $$
in suitably weighted $ L^2 $ or $ L^p $ spaces. The analysis of the resolvent near the thresholds of the spectrum is a natural question in itself, but we shall also discuss  applications to the time decay of wave and Schr\"odinger equations.

We first consider resolvent estimates. The study of the limiting absorption principle, namely the behaviour of (powers of) the resolvent of self-adjoint operators as the spectral parameter approaches the absolutely continuous spectrum is a basic problem in scattering theory and there is a huge literature on this topic which we can not review here.  For the operators considered in this paper (and more general Schr\"odinger operators), the
 analysis of $ (-\Delta_G-z)^{-n} $ is rather well known  as long as $ \mbox{Re}(z) $ remains away from $ 0 $;   by the results of \cite{JeMoPe,AJensen} (and those of \cite{KochTataru} to ensure that $ - \Delta_G $ has no embedded eigenvalues in its (absolutely continuous) spectrum $ [0,\infty) $), we know that, for any $ I \Subset (0,\infty) $ and $ n \geq 1 $, the limits $ \lim_{\epsilon \rightarrow 0^{\pm}} (- \Delta_G- \lambda - i \epsilon)^{-n} $ exist as bounded operators between dual weighted $ L^2 $ spaces, provided 
 that $ \lambda \in I $. 
 The asymptotics as $ \lambda \rightarrow + \infty $ have also been widely studied in various contexts, perhaps more for the resolvent itself than for its powers, but this is not a serious restriction since, in the high energy or semiclassical regime, one can get estimates for powers in terms of estimates of the resolvent (see \cite{IKRem,AJensen} and Subsection \ref{resolvanteverspuissances} below): basically  $ ||(-\Delta_G - \lambda - i 0)^{-n}||^{\prime} $ grows as $ || (- \Delta_G - \lambda -i0)^{-1} ||^n $, if $ || \cdot || $ and $ || \cdot ||^{\prime} $ are operator norms between suitable weighted $ L^2 $ spaces. In this regime, the asymptotics depend crucially on whether the geodesic flow is non trapping, namely if all geodesics escape to infinity as time goes to infinity, or trapping: see \cite{RoTa,Wang,GeMa,RoENS,VaZw} for the non trapping case, \cite{Naka1,NoZw} for weak trapping, and \cite{Burq1,Burq2,stCardosoVodev} in the general case, ie without condition on the geodesic flow.
 
  The situation is definitely different as $ \mbox{Re}(z) \rightarrow 0 $. At first, we note that the geodesic flow plays no role in this non semiclassical regime. More importantly, there is no hope to deduce bounds on powers of the resolvent from bounds on the resolvent $ (-\Delta_G - z)^{-1}  $  as above. We know indeed that $ (-\Delta_G - z)^{-1} $ remains bounded for $ z $ close to $ 0 $ (see \cite{Bouc1,BoHa2,GuHa} for the long range metric case) but, as we shall see below, its powers start to blow up as $ z \rightarrow 0 $ if $n$ is large enough (essentially $ n > d/2 $).
This can be seen on the example of the flat Laplacian on $ \Ra^3 $ whose kernel of the resolvent reads
$$ G_z (x,y) = \frac{1}{4 \pi} \frac{e^{i z^{1/2}|x-y|}}{|x-y|} , \qquad \mbox{Im}(z^{1/2}) > 0 . $$ 
Indeed, since $ (- \Delta -z)^{-2} = \frac{d}{dz} (-\Delta - z)^{-1} $, we  see  that for $ z \in \Ca \setminus \Ra $,
\begin{eqnarray}
 || (- \Delta - z)^{-1} ||_{L^{6/5} \rightarrow L^6} \lesssim 1, \qquad ||(-\Delta-z)^{-2}||_{L^1 \rightarrow L^{\infty}} \approx |z|^{-1/2}, \label{estimeeslibres}
\end{eqnarray}
where the first estimate follows from the Hardy-Littlewood-Sobolev inequality. Of course, such $ L^p \rightarrow L^{p^{\prime}} $ estimates imply weighted $ L^2 $ estimates using, in the present case, the boundedness of
$$ \scal{x}^{-1- \varepsilon} : L^2 (\Ra^3) \rightarrow L^{6/5} (\Ra^3),  \qquad \scal{x}^{- \frac{3}{2}- \varepsilon} : L^2 (\Ra^3)  \rightarrow L^1 (\Ra^3) ,  $$ 
and their adjoints (for any $ \varepsilon > 0 $). 
 
 The literature on powers of the resolvent near the $ 0 $ energy is rather lacunary in the long range case. Actually, this topic seems to have been studied for Schr\"odinger operators $ - \Delta + V $ only, in \cite{Naka2} for $V$ of definite sign and in \cite{FoSk} for $V$ sufficiently negative at infinity (see also \cite{Yafa} in the radial case). We note that, for such potentials, the resolvent behaves differently to the free resolvent in that its powers are uniformly bounded as $ \mbox{Re}(z) \rightarrow 0^+ $, unlike (\ref{estimeeslibres}). Our first purpose is to show that, for variable coefficients metrics, we get the same kind of estimates as in the free case.
 
To state our results, we introduce the notation
\begin{eqnarray}
 \bar{r} (d) = \mbox{the largest integer strictly smaller than } d/2  . \label{notationregularitemaximale}
\end{eqnarray} 
In other words
$ \bar{r}(d) $ is the integer part $ [d/2] $ of $d/2$ if $ d $ is odd, and $ \frac{d}{2}-1 $ if $d$ is even. 
The notation $r$ refers to the fact that it will be interpreted as some regularity index  further on.
Let us remark  that, in all cases, $ \bar{r}(d) \geq 1 $. We also introduce the conjugate Lebesgue exponents
\begin{eqnarray}
p (n) = \frac{2d}{d+2n} , \qquad q (n) = \frac{2d}{d-2n} , \qquad \mbox{for} \ \ 1 \leq n \leq \bar{r}(d) ,\label{notationLebesgue}
\end{eqnarray}
which belong to $ (1,\infty) $ since $ n < d/2 $ by definition of $ \bar{r}(d) $. Finally, we denote by $A$ the (self-adjoint realization of the) generator of $ L^2$ dilations, namely
\begin{eqnarray}
 A  =  \frac{x \cdot \nabla}{i} + \frac{d}{2} . \label{generateur}
\end{eqnarray}
Our main result is the following. 
\begin{theo} \label{estimationprincipaleLp} Fix $ d \geq 3 $. There exists $ \kappa > 0 $ and $ C > 0 $ such that, for $ 1 \leq n \leq \bar{r}(d) $,
\begin{eqnarray}
 \big| \big| (\kappa A + i)^{-n} (- \Delta_G - z)^{-n} (\kappa A - i)^{-n} \big| \big|_{L^{p(n)} \rightarrow L^{q(n)}} \leq C, \qquad |\emph{Re}(z)|<1,  
 \nonumber
\end{eqnarray}
and, for $  n=N:= \bar{r}(d) + 1   $, 
\begin{enumerate}
\item{if $d$ is odd, 
\begin{eqnarray}
\big| \big| (\kappa A + i)^{-N} (- \Delta_G - z)^{-N} (\kappa A - i)^{-N} \big| \big|_{L^1 \rightarrow L^{\infty}} \leq C |\emph{Re}(z)|^{-1/2} , 
\nonumber
\end{eqnarray}}
\item{if $d$ is even, then for all $ q > 2d $ there exists $ C_q > 0 $ such that,
\begin{eqnarray}
\big| \big| (\kappa A + i)^{-N} (- \Delta_G - z)^{-N} (\kappa A + i)^{-N}  \big| \big|_{L^{q/(q-1)} \rightarrow L^{q}} \leq C_{q} |\emph{Re}(z)|^{-\frac{2d}{q}},
\nonumber  
\end{eqnarray}}
\end{enumerate}
both for $ 0<|\emph{Re}(z)|<1  $.
\end{theo}

Up to the weights $ (\kappa A \pm i)^{-n} $, this theorem generalizes the estimates (\ref{estimeeslibres}) to long range metrics in all dimensions greater than 2. In odd dimensions, our result is sharp from the point of view of the singularity at $z=0$, as shown by (\ref{estimeeslibres}).
We also point out that for small long range perturbations of the flat Laplacian (when $ G - I $ is small everywhere on $ \Ra^d $, not only at infinity as imposed by (\ref{cometriquelongueportee})), such estimates actually hold for all $z$, ie also for large ones, and  are scale invariant (see Subsection \ref{subsectionsmall}).

The $ L^p \rightarrow L^{p^{\prime}} $ estimates of Theorem \ref{estimationprincipaleLp} can be turned into the following weighted $ L^2 \rightarrow L^2 $ estimates.
In fact, as one can see from Subsection \ref{nonsmallsubsection} below, Theorem \ref{estimationprincipaleLp} and the following one are equivalent.
\begin{theo} \label{estimationprincipale} Let $ N = \bar{r}(d) + 1 $. For all $ 1 \leq n \leq N  $ and $ \nu > 2n $, we have: for $ n \leq N-1 $,
\begin{eqnarray}
 \big| \big| \scal{x}^{-\nu} (- \Delta_G - z)^{-n} \scal{x}^{-\nu} \big| \big|_{L^2 \rightarrow L^2} \leq C, \qquad |\emph{Re}(z)|<1,  \label{resolventalone}
\end{eqnarray}
and, for $  n= N  $, 
\begin{enumerate}
\item{if $d$ is odd, 
\begin{eqnarray}
\big| \big| \scal{x}^{-\nu} (- \Delta_G - z)^{-N} \scal{x}^{-\nu} \big| \big|_{L^2 \rightarrow L^2} \leq C |\emph{Re}(z)|^{-1/2} ,  \nonumber
\end{eqnarray}}
\item{if $d$ is even, then for all $ \epsilon > 0 $,
\begin{eqnarray}
\big| \big| \scal{x}^{-\nu} (- \Delta_G - z)^{-N} \scal{x}^{-\nu}  \big| \big|_{L^2 \rightarrow L^2} \leq C_{\epsilon} |\emph{Re}(z)|^{-\epsilon}, \nonumber
\end{eqnarray}}
\end{enumerate}
both for $ 0<|\emph{Re}(z)|<1  $.
\end{theo}
We implicitly assume that $ \mbox{Im}(z) \ne 0 $ in these estimates but, since one knows that pointwise limits exist as $  \pm \mbox{Im}(z) \rightarrow 0^+ $ (if $ \mbox{Re}(z) \ne 0 $), this is not a restriction.


\bigskip

 We now study the applications to the local energy decay. Let us consider
\begin{eqnarray}
i \partial_t u + \Delta_G u &= & 0, \qquad u_{|t=0} = u_0 , \label{Schrodinger} \\
\partial_t^2 w - \Delta_G w +m^2 w & = & 0, \qquad w_{|t=0} = f , \ \ \partial_t w_{|t=0} = g, \label{KleinWave}
\end{eqnarray}
which are respectively the Schr\"odinger equation and Klein-Gordon equation ($m > 0$) or wave equation ($ m = 0 $). We are interested in the decay as $ t \rightarrow \infty $ of
\begin{eqnarray}
 || \chi u (t) ||_{L^2}, \qquad || \chi \partial_t w (t) ||_{L^2} + || \chi w(t) ||_{H^1} , \label{energies}
\end{eqnarray}
for some spatial localization $ \chi $, typically $ \chi \in C_0^{\infty} $ or more generally $ \chi (x) = \scal{x}^{-\nu} $ for some $ \nu > 0  $. Using estimates on the resolvent alone (ie those for $n=1$ in (\ref{resolventalone})), it is well known that one can recover $ L^2 (\Ra,dt) $ estimates for (\ref{energies}) which is a weak form of time decay (see for instance \cite{RodTao,Tataru,MeTa,BoHa1,VaWu} in contexts close to ours). Proving quantitative decay rates requires more information, for instance estimates on powers of the resolvent as we recall now.

The flows of the equations (\ref{Schrodinger}) and (\ref{KleinWave}) are functions of $ \Delta_G $ namely,
$$ u(t) = e^{it\Delta_G} u_0, \qquad w (t) = \cos \left( t \sqrt{m^2-\Delta_G} \right) f + \frac{\sin \left( t \sqrt{m^2-\Delta_G} \right) }{ \sqrt{m^2-\Delta_G}}g . $$
If we denote by $ (E_{\lambda})_{\lambda \in \Ra} $ the family of spectral projections associated to $ - \Delta_G $, the latter reads, for instance for the Schr\"odinger equation,
\begin{eqnarray}
 e^{it \Delta_G} = \int e^{-it\lambda} d E_{\lambda} , \label{FourierSchrodinger}
\end{eqnarray}
where the spectral measure $ d E_{\lambda} $ can be recovered from the resolvent by the following Stone's formula (see \cite{ReedSimon} for a proof and Lemma \ref{Stones} below for a precise statement)
$$ d E_{\lambda} = \lim_{\epsilon \downarrow 0} \frac{1}{2 i \pi} \left(  (-\Delta_G-\lambda - i \epsilon)^{-1} -  (-\Delta_G-\lambda + i \epsilon)^{-1} \right) d \lambda . $$
Thus, by using the Stone formula in (\ref{FourierSchrodinger}) and integrating by part, one expects to recover time decay for (\ref{energies}) from the smoothness of the resolvent with respect to $ \lambda $, that is from the integrability in $ \lambda $ of $ \lim_{\epsilon \rightarrow 0^+} (-\Delta_G- \lambda \pm i \epsilon)^{-n} $.

This is of course a rough formal description but this approach is well known and can be made rigorous (see Section \ref{referenceStone}), provided that the corresponding integrals are convergent with respect to $ \lambda $. Its justification requires two types of estimates: high frequency estimates ($ \lambda \rightarrow \infty $) and low frequency estimates ($ \lambda \rightarrow 0 $). As we recalled above, the high frequency estimates are often a delicate question, especially when there are trapped geodesics. But independently of this question, whatever the classical dynamic looks like and whatever methods are used to study the resolvent ({\it e.g.} resonances theory or Mourre theory), one also needs to deal with the low frequencies. For the local energy decay, the latter has been treated in the litterature for  fast decaying perturbations (from \cite{LaPh,Vain,Burq1} for compactly supported perturbations to \cite{Wang2,SSS} in the short range case, with radial assumptions in \cite{SSS}). For long range metrics, there are either conditionnal results (\cite{Christianson} which assume that the resolvent can be continued accross the absolutely continuous spectrum near $ 0 $) or spectrally localized estimates (\cite{CardosoVodev2,Vodev} where the evolution $ e^{it (-\Delta_G)^{\alpha}} $ is replaced by $ e^{it(-\Delta_G)^{\alpha}} \psi (\Delta_G) $ with $ \psi \equiv 0 $ near $0$). We should also mention the recent results for asymptotically flat space times   \cite{Tataru2} where a (sharp) pointwise energy decay is obtained for long range perturbations, which are radial up to short range terms (see also \cite{DSS} for a purely radial situation).


Using mainly Theorem \ref{estimationprincipale} and the known results on high energy estimates \cite{stCardosoVodev}, we  obtain the following general result.
\begin{theo} \label{theoremeexponentiel} Let $d\geq 3 $ and assume (\ref{cometriquelongueportee}). For all positive real numbers $s>0 $, $ \nu > 0 $ there exists $ C > 0 $ such that, for the Schr\"odinger equation (\ref{Schrodinger}),
$$ \big| \big| \scal{x}^{-\nu} u (t) \big| \big|_{L^2} \leq C \left( 1 + \log \scal{t} \right)^{-s} || \scal{x}^{\nu} u_0 ||_{H^s} , $$
and, for the wave and Klein-Gordon equations (\ref{KleinWave}),
$$ \big| \big| \scal{x}^{-\nu} \partial_t w (t) \big| \big|_{L^2} + \big| \big| \scal{x}^{-\nu} w (t) \big| \big|_{H^1} \leq C \left( 1 + \log \scal{t} \right)^{-s} \left( || \scal{x}^{\nu} f ||_{H^{1+s}} + || \scal{x}^{\nu} g ||_{H^s} \right) . $$
\end{theo}
We emphasize that the main novelty in this result is that no spectral cutoff is needed on the initial data and that the metric $G$ is long range. Furthermore, we don't use any spherical symmetry. The time decay is very weak, but on the other hand there is no assumption on the geodesic flow. This is an analogue of a result of Burq \cite{Burq1}, obtained initially for compactly supported perturbations of the Laplacian and then generalized to long range perturbations but with spectrally localized initial data in \cite{CardosoVodev2}.

Under the non trapping condition, we obtain the following stronger decay which we shall prove for the Schr\"odinger equation only. Here we use the notation (\ref{notationregularitemaximale}).
\begin{theo} \label{theoremenontrapping} Let $d\geq 3 $, assume (\ref{cometriquelongueportee}) and that the geodesic flow is non trapping. For all real numbers $ \nu > 2 ( \bar{r}(d) + 1 ) $ and $ 0 \leq s < \nu $, there exists $ C > 0 $ such that, for the Schr\"odinger equation (\ref{Schrodinger}),
\begin{eqnarray}
 \big| \big| \scal{x}^{-\nu} u (t) \big| \big|_{L^2} \leq C \scal{t}^{-\bar{r}(d)} || \scal{x}^{\nu} u_0 ||_{H^{-s}} . \label{effetregularisant}
\end{eqnarray}
\end{theo}
Notice that, in addition to the time decay, (\ref{effetregularisant}) also means that  we  have a smoothing effect (as is naturally expected for the Schr\"odinger equation with a non trapping metric).
 
 As in Theorem \ref{theoremeexponentiel}, the main point in Theorem \ref{theoremenontrapping} is again  the (non radial) long range assumption and the absence of spectral localization on the initial data. Besides we note that if one avoids the low frequencies, ie replaces $ u (t) $ by $  \Phi (\Delta_G) u (t) $ with $ \Phi \equiv 0 $ near $ 0 $ and smooth, one can show that $  \scal{x}^{-\nu} \Phi (\Delta_G) u (t) $ decays as $ t^{\epsilon - \nu} $, ie with time decay rate growing with the spatial decay rate. Thus we have a fast decay in time if $ \scal{x}^{-\nu} $ is replaced by a Schwartz function. This illustrates the fact that, in the non trapping case, the time decay is governed by the low frequency part of the spectrum. From the free case, we also know that this decay cannot be more than $ \scal{t}^{-d/2} $ and we note that, in odd dimension,  $ \scal{t}^{-d/2} = \scal{t}^{-\bar{r}(d)-\frac{1}{2}} $.

 \medskip
 
\noindent {\bf Remark.} We comment that, in principle, our method would show  for the wave and Klein-Gordon equations with non trapping metrics that
$$ \big| \big| \scal{x}^{-\nu} \partial_t w (t) \big| \big|_{L^2} + \big| \big| \scal{x}^{-\nu} w (t) \big| \big|_{H^1} \leq C \scal{t}^{-\bar{r}(d)} \left( || \scal{x}^{\nu} f ||_{H^{1}} + || \scal{x}^{\nu} g ||_{L^2} \right) . $$
We leave this  as a remark since its proof would require bounds of the form, for $ \varphi \in C_0^{\infty}(0,+\infty) $,
$$  || \scal{x}^{-\nu} e^{i t (m^2 - \Delta_G)^{1/2}} \varphi (-h^2\Delta_G) \scal{x}^{-\nu} || \leq C_{\epsilon} \scal{t}^{\epsilon-\nu} , \qquad t \in \Ra , \ h \in (0,1] , $$
for all $ \epsilon > 0 $ (see the proof of Theorem \ref{theoremenontrapping} in Subsection \ref{decroissancetemporelle}). The point in the latter estimate is that it is uniform with respect to $h$ and cannot be clearly deduced from resolvent estimates only. It could however certainly be obtained as the similar one  for the Schr\"odinger equation proved by Wang \cite{Wang} using the Isozaki-Kitada parametrix. 
 




\section{Model operators} \label{section2}
\setcounter{equation}{0}
In this section, we introduce a class of second order differential operator which are small perturbations of the flat Laplacian. They will serve as models at infinity, in the sense  that the Laplacians (\ref{LaplacienG}) will be unitarily equivalent to compactly supported perturbations of such operators (see Subsection \ref{nonsmallsubsection}).
These model operators are of the form
\begin{eqnarray}
 P = a_{jk}(x) D_j D_k +  b_k (x) D_k , \label{formedePdetail}
\end{eqnarray}
where we use the summation convention and $ D_k = i^{-1} \partial/\partial x_k $. We also assume that
\begin{eqnarray}
P \ \mbox{is formally symmetric on} \  L^2(\Ra^d) \ \mbox{with respect to the Lebesgue measure}. \label{conditionsymetrie}
\end{eqnarray}
By formally symmetric, we mean symmetric when tested again functions of the Schwartz space $ {\mathcal S} $. The coefficients will be chosen in the following spaces. 
For integers such that
$$ N \geq 0, \qquad 0 \leq o \leq 1 , \qquad r \geq 0 \qquad \mbox{and} \qquad o + r \leq \bar{r}(d) , $$
where $ \bar{r}(d) $ is defined by (\ref{notationregularitemaximale}), we introduce the norm
$$ || a ||_{o,r,N} = \sum_{n \leq N, \atop |\alpha | \leq r } || \partial^{\alpha} (x \cdot \nabla)^n a ||_{L^{\frac{d}{o+|\alpha|}}} , $$
and define
$$ S^{o,r,N} = \{ a \in C^{\infty}_b (\Ra^d) \ | \ ||a||_{o,r,N} < \infty \} . $$

In most statements of Sections \ref{section2} and \ref{section3}, we shall also assume that $ P $ is a small perturbation of $ - \Delta $ in the sense that, $ a_{jk} - \delta_{jk}$ ($ \delta_{jk} = $ Kronecker symbols) and $ b_k $ will be small in appropriate spaces, explicitly given in each proposition.  A first example of such a statement is the following.

\begin{prop}  For all $ a_{jk},b_k \in C^{\infty}_b$ such that
$$ \sum_{j,k} ||a_{jk} - \delta_{jk} ||_{L^{\infty}} $$ is small enough, then $ P : {\mathcal S} \rightarrow L^2 $ has a bounded closure $ \overline{P} : H^2 \rightarrow L^2 $ and $ \overline{P} $ is self-adjoint on $ L^2  $ with domain $ H^2 $.
\end{prop}

This proof of this proposition is completly standard  and the existence of a self-adjoint realization holds under much more general assumptions. The smallness of $ \epsilon $ is only used to ensure that the operator is uniformly elliptic. The assumption that the coefficients belong to $ C^{\infty}_b $ guarantees the existence of the closure of $ P $ to  $ H^2 $ and the fact that the domain of $ \overline{P}^* $ is also $ H^2 $ by elliptic regularity. This proposition has to be considered as an algebraic preliminary, which is convenient for it gives explicitly the domain of $ \overline{P} $. But, as far as estimates are concerned, all the bounds obtained below will be given  in terms of $ S^{o,r,N} $  norms only, so the condition $ a_{jk},b_k \in C^{\infty}_b $ is essentially irrelevant.

\medskip

Most of our estimates rely on the following elementary proposition.

\begin{prop} \label{Sobolevhomogene1} There exists $ C > 0$ depending only on the dimension $d$ such that
\begin{eqnarray}
 ||a \partial_j \partial_k u ||_{L^2} \leq C || a ||_{0,1,0} || \Delta u ||_{L^2}, \label{1explicite} \\
 || b \partial_k u ||_{L^2} \leq C || b ||_{1,0,0} ||\Delta u||_{L^2} , \label{2explicite}
 \end{eqnarray}
 and
 \begin{eqnarray}
 ||a \partial_j \partial_k u ||_{H^{-1}} \leq C || a ||_{0,1,0} || \nabla u ||_{L^2}, \label{1explicitebis} \\
 || b \partial_k u ||_{H^{-1}} \leq C || b ||_{1,0,0} || \nabla u||_{L^2} , \label{2explicitebis}
 \end{eqnarray}
 for all $ u \in {\mathcal S} $, all $ a \in S^{0,1,0} $ and $  b \in S^{1,0,0}$.
\end{prop}

\noindent {\it Proof.} We consider first (\ref{1explicite}) and (\ref{2explicite}). By a standard limiting argument we may assume that the Fourier transform of $u$ vanishes near $ 0 $. Then
$$ a (x) \partial_j \partial_k u =  a (x) \frac{\partial_j \partial_k}{|D^2|} |D|^2 u , $$
so (\ref{1explicite}) follows from the bound $ ||a||_{0,1,0} \leq ||a||_{L^{\infty}} $ and the $ L^2  $ boundedness of  $ \partial_j \partial_k / |D|^2 $. To prove (\ref{2explicite}), we use the H\"older inequality
\begin{eqnarray}
 || \varphi \psi ||_{L^2} \leq  || \varphi ||_{L^d} || \psi ||_{L^{2^*}} , \label{Holderautoadjoint}
\end{eqnarray}
and the Sobolev inequality
\begin{eqnarray}
|| \psi ||_{L^{2^* }} \lesssim || \nabla \psi ||_{L^2} . \label{Sobolevautoadjoint}
\end{eqnarray}
Indeed, by writing
$$  b_k(x) \partial_k u = b_k (x)  \frac{1}{|D|} \frac{\partial_k}{|D|} |D|^2 u , $$
 we get the result since  on one hand $ \partial_k / |D| $ is bounded on $ L^2 $  and on the other hand (\ref{Holderautoadjoint})  and (\ref{Sobolevautoadjoint})  yield
$$ || b_k |D|^{-1} \psi ||_{L^2} \leq ||b_k||_{L^d} || |D|^{-1} \psi ||_{L^{2^*}} \lesssim ||b_k ||_{L^d} || \psi ||_{L^2} , $$
for all $ \psi \in {\mathcal S} $ with Fourier transform vanishing near $ 0 $. We now prove (\ref{2explicitebis}). The latter simply follows from the fact that
$$ | (\psi,b \partial_k u) | \leq ||\overline{b} \psi ||_{L^2} || \partial_k u ||_{L^2} \lesssim ||b||_{L^d} || \nabla \psi ||_{L^2} || \nabla u ||_{L^2} , $$ 
using again (\ref{Holderautoadjoint}) and (\ref{Sobolevautoadjoint}). Finally for (\ref{1explicitebis}), we write 
$$ a (x) \partial_j \partial_k = \partial_j (a(x) \partial_k) - (\partial_j a)(x) \partial_k , $$
for which the contribution of the first term follows by an integration by part, and the contribution of the second term follows from (\ref{2explicitebis}) since $ \partial_j a \in S^{1,0,0} $. \finpreuve
\bigskip

A first consequence is the following.

\begin{prop} \label{borneautoadjoint} For all $ P $  satisfying (\ref{conditionsymetrie})  and such that
$$ \sum_{j,k} ||a_{jk} - \delta_{jk} ||_{0,1,0} + ||b_k ||_{1,0,0} $$ is small enough, we have 
\begin{eqnarray}
 \frac{1}{2} || \nabla u ||^2_{L^2} \leq ( \overline{P} u , u ) \leq 2 || \nabla u ||_{L^2}^2 , \label{faible} \\
 \frac{1}{2} || \Delta u ||_{L^2} \leq || \overline{ P } u ||_{L^2}  \leq 2 || \Delta u ||_{L^2}^2 . \label{fort}
\end{eqnarray}
for all $ u \in H^2 $. In particular,
\begin{eqnarray}
\overline{P} \geq 0 . \label{semibornePbar}
\end{eqnarray}
\end{prop}

\noindent {\it Proof.} It suffices to prove the result when $ u \in {\mathcal S} $. To prove (\ref{faible}), we write first
$$ P = D_j a_{jk}(x)  D_k +  \left( b_k (x) - i \partial_j a_{jk}(x) \right) D_k . $$
Then
$$ (P u , u )= || \nabla u ||^2 + \big( (a_{jk}-\delta_{jk})D_k u , D_j u \big) + \big(  (\overline{b_k} + i \partial_j \overline{a_{jk}} ) u , D_k u \big) , $$
so that
$$ (Pu,u) - || \nabla u ||^2 \leq \sum_{j,k} ||a_{jk} - \delta_{jk} ||_{L^{\infty}} || \nabla u ||_{L^2}^2 + 
|| b_k - i \partial_j a_{jk} ||_{L^d} || u ||_{L^{2^*}} || \nabla u ||_{L^2} , $$
using the H\"older inequality (\ref{Holderautoadjoint}).
We conclude with the Sobolev inequality (\ref{Sobolevautoadjoint}). The estimate (\ref{fort}) follows simply from the fact that
 
$$ P u =  - \Delta u +  (a_{jk}(x)-\delta_{jk})D_j D_k u + b_k (x)D_ku , $$
and Proposition \ref{Sobolevhomogene1}. \finpreuve

\bigskip

We next recall the definition and some elementary properties of $A$ the generator of $ L^2 $ dilations
$$ e^{i \tau A} \varphi = e^{\frac{d}{2} \tau} \varphi \big( e^{\tau} x \big) , $$
which is given by (\ref{generateur}). We have the identities
\begin{eqnarray}
|| e^{i \tau A} \varphi ||_{L^p} & = & e^{\tau \left( \frac{d}{2}-\frac{d}{p} \right)} ||\varphi||_{L^p}, \qquad p \in [1,\infty], \label{pourlescalingfinal} \\
|| \partial^{\alpha} e^{i \tau A} \varphi ||_{L^2} & = & e^{\tau |\alpha|} || \partial^{\alpha} \varphi ||_{L^2} ,
\end{eqnarray}
which are convenient to prove the $ L^p \rightarrow L^p $ or $ H^s \rightarrow H^s $ boundedness of 
\begin{eqnarray}
(\kappa A-\zeta)^{-1} = \frac{1}{i} \int_0^{\pm \infty} e^{-i \tau \zeta} e^{i \tau \kappa A} d \tau, \qquad \pm \mbox{Im}(\zeta) < 0 , \label{explicitee}
\end{eqnarray}
for suitable parameters $ \kappa , \zeta,s $ and $p$. For instance, if $ s \geq 0 $ is an integer, we have the useful estimate
\begin{eqnarray}
|| e^{i\tau \kappa A} \varphi ||_{H^s} \leq C \left( 1 +e^{\tau s \kappa} \right) || \varphi ||_{H^s}, \qquad \varphi \in {\mathcal S} .  \label{controleexponentiel}
\end{eqnarray}

\begin{lemm} \label{regulariteSobolevdilatation} There exists $ C>0 $ such that, for all $ r \in \Ra $, $ \kappa > 0 $ and $ \zeta \in \Ca \setminus \Ra $ such that $ |\emph{Im}(\zeta)| > \kappa |r| $, one has
$$ \big| \big|  (\kappa A-\zeta)^{-1} \varphi \big| \big|_{H^r} \leq \frac{C}{|\emph{Im}(\zeta)|- \kappa |r|} ||\varphi||_{H^r}  , $$
for all $ \varphi \in H^r  \cap L^2 $.
\end{lemm}
 
\noindent {\it Proof.} We may assume that $ r \geq 0$, otherwise we consider the adjoint. We consider the norm
$$ || u ||_{H^r} = || u ||_2 + \big| \big| |D|^r u \big| \big|_2  , $$
and recall that
\begin{eqnarray}
 |D|^r ( \kappa A-\zeta)^{-1} = (\kappa A-\zeta-i \kappa r)^{-1} |D|^r  , \label{passagederesolvante}
\end{eqnarray}
which follows from (\ref{explicitee}) (see \cite{Bouc1}).
This formula and the self-adjointness of $A$ give
$$ \big| \big| |D|^r (\kappa A- \zeta	)^{-1} \varphi \big| \big|_{L^2} \leq \frac{1}{|\mbox{Im}(\zeta	)|- \kappa r} \big| \big| |D|^r \varphi \big| \big|_{L^2} , $$
so the result follows using $ \big| \big|  (\kappa A-\zeta)^{-1} \varphi \big| \big|_{L^2} \leq |\mbox{Im}(\zeta	)|^{-1} || \varphi||_{L^2}  $ and $ |\mbox{Im}(\zeta	)|^{-1} \leq (|\mbox{Im}(\zeta)|-\kappa r)^{-1} $.
 \finpreuve 

\bigskip

We now consider commutators with $A$. We recall the following standard notation, 
$$ \mbox{ad}^0_A P = P, \qquad \mbox{ad}_A P = [P,A], \qquad \mbox{ad}^n_A P = \big[ \mbox{ad}_A^{n-1} P , A \big] . $$
Here the commutators are defined in the sense of differential operators acting on Schwartz functions. One easily checks that
\begin{eqnarray}
 i^n \mbox{ad}_A^n P = a_{jk}^{(n)}(x) D_j D_k +  b_k^{(n)} (x) D_k , \label{niemecommutateur}
\end{eqnarray}
with
\begin{eqnarray}
 a_{jk}^{(n)}= \big(2 - x \cdot \nabla \big)^n a_{jk}, \qquad  b_k^{(n)} = \big( 1 - x \cdot \nabla \big)^n b_k . \label{formedescoefficients}
\end{eqnarray}

\begin{prop} For all $ n \geq 0 $, there exists $ C_n $ such that
\begin{eqnarray}
 || \emph{ad}_A^n P u ||_{L^2} \leq C_n \left( \sum_{jk} || a_{jk}  ||_{0,1,n} + ||b_k ||_{1,0,n} \right) ||\Delta u||_{L^2} , \label{Sobolevhomogenecomm}
\end{eqnarray}
for all $ u \in {\mathcal S} $, $ a_{jk} \in S^{0,1,n} $ and $ b_k \in S^{1,0,n} $. In particular,
\begin{enumerate}
\item{If $ \sum_{j,k} ||a_{jk} - \delta_{jk} ||_{0,1,0} + ||b_k ||_{1,0,0} $
is small enough, then
\begin{eqnarray}
 || \emph{ad}_A^n P u ||_{L^2} \leq C_n \left( \sum_{jk} || a_{jk}  ||_{0,1,n} + ||b_k ||_{1,0,n} \right)  || P u ||_{L^2}  , \label{Pborne}
\end{eqnarray}}
\item{ if $ \sum_{j,k} ||a_{jk} - \delta_{jk} ||_{0,1,1} + ||b_k ||_{1,0,1} $
is small enough, then
\begin{eqnarray}
 \big( u, i [P,A] u \big) \geq \frac{1}{2} (u,Pu) , \label{commutateurexplicitenouveau}
\end{eqnarray}}
\end{enumerate}
for all $ u \in {\mathcal S} $.
\end{prop}

\noindent {\it Proof.} The estimate (\ref{Sobolevhomogenecomm}) follows from Proposition \ref{Sobolevhomogene1} and (\ref{formedescoefficients}). If $ \sum_{j,k} ||a_{jk} - \delta_{jk} ||_{0,1,0} + ||b_k ||_{1,0,0}  $ is small then we may replace $ || \Delta u ||_{L^2} $ by $ ||Pu ||_{L^2} $ in (\ref{Sobolevhomogenecomm})  using (\ref{fort}), which proves (\ref{Pborne}). To prove (\ref{commutateurexplicitenouveau}), we proceed  similarly to (\ref{faible}) to show that $ \big( u, i [P,A] u \big) \geq || \nabla u ||_{L^2}^2 $ (note that $ i [P,A]$ is close to $ - 2\Delta $) and use (\ref{faible}) to conclude. \finpreuve

\bigskip

 The estimate (\ref{commutateurexplicitenouveau}) is a positive commutator estimate which holds uniformly for all $ a_{jk} ,b_k $ in bounded subsets of $ S^{0,1,n} $ and $ S^{1,0,n}$ respectively and satisfying the smallness condition of item 2.  In the same spirit, the estimate (\ref{Pborne}) means that $ \mbox{ad}_A^n P $ is relatively bounded with respect to $ P $ with a fairly explicit dependence on the coefficients $a_{jk}$, $b_k$. In the next proposition, we derive some useful related estimates.

\begin{prop} \label{fermeturecommutateursprop} For all $n$, $ \emph{ad}^n_A P : {\mathcal S} \rightarrow L^2 $ has a bounded closure 
$$ \overline{\emph{ad}^n_A P} : H^2 \rightarrow L^2 , $$
provided that $ a_{jk} \in S^{0,1,n} $, $ b_k \in S^{1,0,n} $. If in addition (\ref{conditionsymetrie}) holds and
\begin{eqnarray}
 \sum_{j,k} ||a_{jk} - \delta_{jk} ||_{0,1,0} + ||b_k ||_{1,0,0} \label{smallrelecture} 
\end{eqnarray}
is small enough, then for all integer $ j \geq 0 $,
\begin{eqnarray}
 \big| \big| \overline{\emph{ad}^n_A  P} (\kappa A+i)^{-j} \big( \overline{P} - z \big)^{-1} \big| \big|_{L^2 \rightarrow L^2} \leq \frac{C_{n,j}}{ (1-2 \kappa)^{j}} \left( \sum_{jk} || a_{jk} ||_{0,1,n} + ||b_k ||_{1,0,n} \right) \frac{\scal{z}}{{\rm dist}(z,[0,\infty))} ,  
 \nonumber \label{estimeespectrale}
\end{eqnarray}
for all $ z \in \Ca \setminus [ 0, \infty ) $ and all $ 0 < \kappa < 1/2 $. In particular, if the coefficients $ a_{kj} ,b_k $ belong to bounded subsets of $ S^{0,1,n} $ and $ S^{1,0,n} $ respectively, and if (\ref{smallrelecture}) is small enough, then $ \overline{\emph{ad}^n_A P} $ is $ \overline{P}$ bounded and
\begin{eqnarray}
 \big| \big| \overline{\emph{ad}^n_A  P}  \big( \overline{P} +1 \big)^{-1} \big| \big|_{L^2 \rightarrow L^2} \lesssim 1 , \label{Pborneuniforme}
\end{eqnarray}
 uniformly with respect to these coefficients.
\end{prop}

\noindent {\it Proof.} The existence of the closure follows from (\ref{Sobolevhomogenecomm}). The prove the estimate, we write
\begin{eqnarray*}
  \overline{\mbox{ad}^n_A  P} (\kappa A+i)^{-j} \big( \overline{P} - z \big)^{-1} & = &   \overline{\mbox{ad}^n_A  P} (1-\Delta)^{-1} \left( (\kappa A+i)^{-j} \big( \overline{P} - z \big)^{-1} \right. \\
  & & \left. - (\kappa A+i-2i \kappa)^{-j} \Delta \big( \overline{P} - z \big)^{-1} \right) .
\end{eqnarray*}
using (\ref{passagederesolvante}).
We conclude by using (\ref{Sobolevhomogenecomm}) and the lower bound in (\ref{fort})  which shows that, for $ k = 0 ,1 $, 
$$ \big| \big| \Delta^k (\overline{P} - z)^{-1} \big| \big|_{L^2 \rightarrow L^2} \lesssim \big| \big| \overline{P}^k (\overline{P} - z)^{-1} \big| \big| \leq \sup_{\lambda \in [ 0 , \infty ) } \left| \frac{\lambda^k}{\lambda - z} \right| \lesssim \frac{\scal{z}^k}{{\rm dist}(z,[0,\infty))}  , $$
the second estimate following from the Spectral Theorem and the fact that $ \mbox{spec}(\overline{P}) \subset [0,\infty) $ by (\ref{semibornePbar}). \finpreuve

\section{Weighted functionnal calculus} \label{section3}
\setcounter{equation}{0}
In this section, we investigate the $ L^2 \rightarrow L^2 $ (and sometimes $ H^{-1} \rightarrow H^1 $) boundedness  of operators of the form
$$ (\kappa A+i)^{-n} \chi \big( \overline{P} \big) (\kappa A+i)^{n} , $$
with $ P $ as in Section \ref{section2} and 
where $ \chi $ may be a bump function in $ C_0^{\infty} $, or corresponds to $ (\overline{P}-z)^{-1} $ ie $ \chi (\alpha) = (\alpha -z)^{-1} $. The expression above is not well defined on $ L^2 $ for it can be applied only to functions in $ \mbox{Dom}(A^n) $, but we shall see that it has a bounded closure. This kind of results is well known, but the additional point we want to stress here is to which extent the norms of these operators are uniform with respect to the coefficients $ a_{jk},b_k $ defining $ P $ in (\ref{formedePdetail}). Since this is a crucial tool in the proof of Theorem \ref{estimationprincipale} (or, more precisely, of Theorem \ref{theoremeprincipalpetit} below), we devote a section to this topic.

\bigskip

The following lemma will be of constant use in the sequel.
\begin{lemm} \label{Lemme1} For all $ n \geq 0 $, all $ u \in H^s $, with $ s \geq 0 $, and all $ \kappa > 0$ such that $ \kappa s < 1 $, there exists a sequence $ \theta_j $ in $ C_0^{\infty} $ such that, for $ k = 0 , \ldots , n $,
$$  (\kappa A+i)^k \theta_j \rightarrow (\kappa A+i)^{k-n} u, \qquad \mbox{in} \ H^s, $$
as $ j \rightarrow \infty $.
\end{lemm}

\noindent {\it Proof.} Let first $ T_j = \frac{1}{i} \int_0^j  e^{-\tau} e^{i \tau \kappa A} dt $.
By (\ref{controleexponentiel}), 
$$ T_j \rightarrow (\kappa A+i)^{-1}, \qquad \mbox{strongly on} \ H^s , $$ and
 using $ \kappa A e^{it\kappa A} = -i \frac{d}{dt} e^{it\kappa A} $,
$$ (\kappa A+i) T_j = I - e^{-j} e^{ij\kappa A} \qquad \mbox{on} \ \mbox{Dom}(A) , $$
where
$$ I - e^{-j} e^{ij\kappa A} \rightarrow I, \qquad \mbox{strongly on} \ H^s . $$
We can now prove the existence of $ \theta_j $ by induction on $n$. The result is clear if $n = 0$. For $ n \geq 1 $, the induction assumption allows to pick $ \varphi_j $ in $ C_0^{\infty} $ such that, for $ k \leq n-1 $,
$$ (\kappa A+i)^k \varphi_j \rightarrow (\kappa A+i)^{k-n+1} u $$ in $ H^s $. We then define
$$  \theta_j = T_j   \varphi_j . $$
Clearly,  $ \theta_j  $ belongs to $ C_0^{\infty} $ since the integral defining $ T_j $ is over a bounded interval.  Furthermore, since $ (\kappa A+i)^k $ commutes with $ T_j $, we have, for $ k \leq n-1 $,
\begin{eqnarray}
 (\kappa A+i)^k \theta_j & = & T_j (\kappa A+i)^k \varphi_j \nonumber \\
  & = & T_j (\kappa A+i)^{k-n+1} u + o (1) \ \longrightarrow (\kappa A+i)^{k-n} u . \nonumber
\end{eqnarray}
Further, for $ k = n $,
\begin{eqnarray*}
 (\kappa A+i)^n \theta_j &= & (\kappa A+i) T_j (\kappa A+i)^{n-1} \varphi_j , \\
& = & (\kappa A+i)^{n-1} \varphi_j - e^{-j} e^{ij\kappa A} (\kappa A+i)^{n-1} \varphi_j \ \longrightarrow u , 
\end{eqnarray*}
and the result follows. \finpreuve

\bigskip

\begin{defi} \label{definitionfermeture} If $ B : L^2 \rightarrow L^2 $ is a bounded operator and $ n \geq 0 $ an integer such that 
$$ (\kappa A+i)^{-n} B (\kappa A+i)^n : {\mathcal S} \rightarrow L^2 ,$$
has a bounded closure to $ L^2 $, we denote by
$$ B_{\kappa,A,n} = \overline{ (\kappa A+i)^{-n} B (\kappa A+i)^n } , $$
this $ L^2 \rightarrow L^2 $ closure.
\end{defi}

\begin{prop} \label{Propositionformelle0} Let $ B $ be a bounded operator such that $ B_{\kappa, A,n} $ exists. Then
\begin{enumerate}
\item{$$ B_{\kappa, A,n} (\kappa A+i)^{-n} = (\kappa A+i)^{-n} B . $$}
\item{If $ C $ is another bounded operator such that $  C_{\kappa, A,n} $ exist then $ (BC)_{\kappa, A,n} $ exists as well and
$$ B_{\kappa, A,n} C_{\kappa, A,n} = (BC)_{\kappa, A,n} . $$}
\end{enumerate} 
\end{prop}
Item {\it 2} gives a rigorous sense to the formally trivial identity
$$ (\kappa A+i)^{-n} B (\kappa A + i)^{n} (\kappa A+i)^{-n} C (\kappa A+i)^n = (\kappa A+i)^{-n} BC (\kappa A+i)^n . $$

\bigskip

\noindent {\it Proof.} To prove that {\it 1} holds when applied on any $ u \in L^2 $, we use Lemma \ref{Lemme1} to pick $ \theta_j \in C_0^{\infty} $ which approaches $ (\kappa A+i)^{-n} u $ and such that $ (\kappa A+i)^n \theta_j $ approaches $u$, both in $ L^2 $.

We prove now {\it 2}. It suffices to show that, for all $ \psi \in {\mathcal S} $,
\begin{eqnarray}
 (\kappa A+i)^{-n} BC (\kappa A+i)^n \psi = B_{\kappa, A,n} C_{\kappa, A,n} \psi . \label{psisuffisant}
\end{eqnarray}
Fix such a $ \psi $ and let
$$ u = C_{\kappa, A,n} \psi = (\kappa A+i)^{-n} C (\kappa A+i)^n \psi \in L^2 . $$
Choose  $ \theta_j $ as in Lemma \ref{Lemme1} so that $ \theta_j \rightarrow u $ and $ (\kappa A+i)^n \theta_j \rightarrow C (\kappa A+i)^n \psi $ in $ L^2 $. Then, on one hand
\begin{eqnarray}
 B_{\kappa, A,n} \theta_j \rightarrow B_{\kappa, A,n} C_{\kappa, A,n} \psi , \label{limite1}
\end{eqnarray}
and on the other hand,
\begin{eqnarray*}
(\kappa A+i)^n B_{\kappa, A,n} \theta_j = B (\kappa A+i)^n \theta_j \rightarrow  BC (\kappa A+i)^n \psi, \qquad \mbox{in} \ L^2 , 
\end{eqnarray*}
so we get
\begin{eqnarray*}
 B_{\kappa, A,n}  \theta_j  \rightarrow  (\kappa A+i)^{-n} BC (\kappa A+i)^n   \psi, \qquad \mbox{in} \ L^2 ,
\end{eqnarray*}
which, together with (\ref{limite1}), implies (\ref{psisuffisant}). \finpreuve

\bigskip

For future purposes, we also record the following straightforward lemma which gives a precise meaning  to the formal expression
$$ (\kappa A+i)^{-n} B (\kappa A + i)^{n} = (\kappa A+i)^{-1} \left( (\kappa A+i)^{1-n} B (\kappa A + i)^{n-1} \right) (\kappa A+i) . $$
\begin{lemm} \label{Lemmeformel1} Let $ B $ be such that $ B_{\kappa, A,n} $ and $ B_{\kappa, A,n-1} $ exist. Then
\begin{eqnarray}
 B_{\kappa, A,n} = \overline{ (\kappa A+i)^{-1} B_{\kappa, A,n-1} (\kappa A+i) } , \label{formel1}
\end{eqnarray}
the right hand side denoting the $ L^2 \rightarrow L^2  $ closure of the corresponding operator defined on $ {\mathcal S} $.  
\end{lemm}

\bigskip 

We shall also need the following result.

\begin{prop} \label{Lemme2} Let $ Q $ be a second order differential operator with smooth coefficients such that $ Q $ and $ [Q,A] $, defined on $ {\mathcal S} $, have bounded closures
$$ \overline{Q}, \ \overline{[Q,A]} : H^2 \rightarrow L^2 . $$
Then, for all $ 0 < \epsilon < 1/2 $ and $ u \in H^2$, we have
$$ (\epsilon A+i)^{-1} \overline{Q} u = \overline{Q} (\epsilon A+i)^{-1} u - \epsilon (\epsilon A+i)^{-1} \overline{[Q,A]} (\epsilon A + i)^{-1} u . $$
\end{prop}

Note  that, by Proposition \ref{fermeturecommutateursprop}, any $ Q $ of the form $ \mbox{ad}_A^j P $ satisfies the assumptions of this proposition.

Recall also Lemma \ref{regulariteSobolevdilatation} which shows that $ (\epsilon A + i)^{-1} $ is bounded on $ H^2 $ so that $ \overline{Q} (\epsilon A+i)^{-1} $ and $\overline{[Q,A]} (\epsilon A + i)^{-1}$ are well defined on $H^2$.

\bigskip

\noindent {\it Proof.} Choose $ \theta_j $ as in Lemma \ref{Lemme1}, such that $ (\epsilon A + i)\theta_j \rightarrow u $ and $ \theta_j \rightarrow (\epsilon A + i)^{-1} u $ in $ H^2 $. Observe that
$$ Q (\epsilon A + i) \theta_j = (\epsilon A + i) Q \theta_j + \epsilon [Q,A] \theta_j , $$
and  apply $ (\epsilon A + i)^{-1} $ to this equality. The result follows by letting $ j \rightarrow \infty $. \finpreuve

\bigskip

Applying Proposition \ref{Lemme2} with $ Q = P $, and applying $ (\overline{P}-z)^{-1} $ to the left of the corresponding identity we get:
\begin{lemm} \label{Lemme3} For all $ 0 < \epsilon < 1/2 $ and $ z \notin [0,\infty) $,
$$ (\epsilon A + i)^{-1} \big( \overline{P} - z \big)^{-1} = \big( \overline{P} - z \big)^{-1} (\epsilon A + i)^{-1} - \epsilon 
\big( \overline{P} - z \big)^{-1} (\epsilon A + i)^{-1} \overline{[P,A]} (\epsilon A + i)^{-1} \big( \overline{P} - z \big)^{-1} ,  $$
as operators from $ H^2 $ to $ L^2 $.
\end{lemm}

The latter lemma is useful to prove the following identity (note that we swap the resolvents of $ \overline{P} $ and $A$).
\begin{prop} \label{Lemme4} For all $ \psi \in C_0^{\infty} $, all $ \kappa > 0$ and all $ z \notin [0,\infty) $, we have
\begin{eqnarray*}
 (\kappa A+i)^{-1} \big( \overline{P} - z \big)^{-1} (\kappa A + i) \psi & = & \big( \overline{P} - z \big)^{-1} \psi +  
  \kappa (\kappa A+i)^{-1} \big( \overline{P} - z \big)^{-1} \overline{[P,A]} \big( \overline{P} - z \big)^{-1} \psi .
\end{eqnarray*} 
\end{prop}

\noindent {\it Proof.} By the Spectral Theorem, we have
$$ i (\epsilon A + i)^{-1} \rightarrow I, \qquad \epsilon \rightarrow 0 , $$
in the strong sense on $ L^2 $ but also in $ H^2 $ by (\ref{passagederesolvante}). On the other hand, one easily checks that
$$  i (\epsilon A + i)^{-1} (\kappa A+i)  =   \frac{i \kappa}{\epsilon} I + \left(1 - \frac{\kappa}{\epsilon}  \right) (\epsilon A + i)^{-1}  , $$
so using Lemma \ref{Lemme3},  
\begin{eqnarray*}
 \big( \overline{P} - z \big)^{-1} (\kappa A + i) i (\epsilon A + i)^{-1} \psi  & = &   (\kappa A + i) i (\epsilon A + i)^{-1}  \big( \overline{P} - z \big)^{-1}  \psi + \\
& & (\kappa - \epsilon) \big( \overline{P} - z \big)^{-1} (\epsilon A + i)^{-1} \overline{[P,A]}  (\epsilon A + i)^{-1} \big( \overline{P} - z \big)^{-1} \psi .
\end{eqnarray*}
Applying $ (\kappa A+i)^{-1} $ to this identity and letting $ \epsilon \rightarrow 0 $, we get the result. \finpreuve
 
\bigskip

\begin{coro} \label{Lemme5} For all $ z \notin [ 0, \infty ) $ and all $ \kappa > 0 $, $ \big( \overline{P} - z \big)^{-1}_{\kappa, A,1} $ exists and is given by
$$ \big( \overline{P} - z \big)^{-1}_{\kappa, A,1} =  \big( \overline{P} - z \big)^{-1}  +  
  \kappa (\kappa A+i)^{-1} \big( \overline{P} - z \big)^{-1} \overline{[P,A]} \big( \overline{P} - z \big)^{-1}  . $$ 
\end{coro}

Notice that we do not need $ \kappa $ to be small since this identity makes sense for operators on $ L^2 $. If we want this result to hold in the sense of operators from $ L^2 $ to $ H^2 $ we have to restrict to $ 0 < \kappa < 1/2 $. 

\bigskip

We next will prove more generally that $ \big( \overline{P} - z \big)^{-1}_{\kappa, A,n} $ exists for any $n$. We will proceed by induction using Lemma \ref{Lemmeformel1}.

\begin{prop} \label{algebrecombinatoire} For all $ z \notin [ 0, \infty ) $, all $ 0 < \kappa < 1/2 $ and all $ n \geq 0 $, $ \big( \overline{P} - z \big)^{-1}_{\kappa, A,n} $ exists and is a linear combination of operators of the form
$$ (\kappa A+i)^{-i_l}   \big( \overline{P} - z \big)^{-1}  \prod_{\nu = 1}^l \left( (\kappa A+i)^{-j_{\nu}} \overline{ \emph{ad}_A^{j_{\nu}} P } (\kappa A+i)^{-k_{\nu}} \big( \overline{P}-z \big)^{-1} \right)  ,
    $$
 the product meaning composition of operators, from the left to the right increasingly in $ \nu $ (it is $I$ if $ l = 0 $), and where 
$$ 0 \leq l \leq n, \qquad 0 \leq  i_l, j_{\nu}, k_{\nu} \leq n.   $$ 
  The coefficients of this combination are non negative powers of $ \kappa $ times complex numbers which are independent of $ \kappa $, $z$ and $ P $.
\end{prop}

\noindent {\it Proof.} We proceed by induction on $n$, the result being trivial if $n=0$. To go from step $n-1$ to $n$, using Proposition \ref{Propositionformelle0} and Lemma \ref{Lemmeformel1},  we have to show that $ B_{\kappa, A,1} $ exist for operators $B$ of the form
$$ (\kappa A+i)^{-k}, \qquad \big( \overline{P} - z \big)^{-1}, \qquad \overline{ \mbox{ad}_A^{j} P } (\kappa A+i)^{-k} \big( \overline{P}-z \big)^{-1} . $$
This is trivial for the first one and follows from Corollary \ref{Lemme5} for the second one. We thus consider the third one, which requires $ \kappa < 1/2 $ to ensure that $ (\kappa A+i)^{-k} $ maps $ H^2 $ in $ H^2 $. By Proposition \ref{Lemme2} (with $ \epsilon = \kappa $), we have
$$ (\kappa A+i)^{-1} \overline{\mbox{ad}_A^j P} = \left( \overline{\mbox{ad}_A^j P} - \kappa (\kappa A+i)^{-1} \overline{\mbox{ad}_A^{j+1} P} \right) (\kappa A+i)^{-1} , $$
which, by Proposition \ref{Propositionformelle0} and Corollary \ref{Lemme5}, shows  that 
$$ \left( \overline{\mbox{ad}_A^j P} (\kappa A+i)^{-k} \big( \overline{P}-z \big)^{-1} \right)_{\kappa, A,1} = \left( \overline{\mbox{ad}_A^j P} - \kappa (\kappa A+i)^{-1} \overline{\mbox{ad}_A^{j+1} P} \right) (\kappa A+i)^{-k} \big( \overline{P}-z \big)^{-1}_{\kappa, A,1} , $$
which is a linear combination of products of operators of the expected form. \finpreuve

\bigskip

We summarize the result obtained so far and derive somes estimates in the following proposition.
\begin{prop} \label{Propositionsansh} There exists $ \epsilon > 0 $ such that, for all integer $ n \geq 0 $ and all $ M > 0 $, there exists $ C > 0 $  such that
  for all coefficients $ a_{jk} \in S^{0,1,n} $, $ b_k \in S^{1,0,n} $ such that
 \begin{enumerate}
 \item{(\ref{conditionsymetrie}) is satisfied}
 \item{$ \sum_{j,k} ||a_{jk} - \delta_{jk}||_{0,1,0} + ||b_k||_{1,0,0} < \epsilon $ }
 \item{$ \sum_{j,k} ||a_{jk} ||_{0,1,n} + ||b_k||_{1,0,n} \leq M $}
\end{enumerate}
and for  all $ z \notin [ 0 , \infty ) $, all $ 0 < \kappa \leq 1/4 $,  we have
\begin{eqnarray}
 \big| \big| (\overline{P}-z)^{-1}_{\kappa, A,n}  \big| \big|_{L^2 \rightarrow L^2} \leq C \left( \frac{\scal{z}}{ {\rm dist}(z,[0,\infty)) } \right)^{n+1 } . \label{sansh0}
\end{eqnarray}
\end{prop}

Notice that we could consider $ 0 < \kappa < 1/2 $, but we restrict to the case $ \kappa \leq 1/4 $ to get a $ \kappa $ independent estimate in (\ref{sansh0}). We would otherwise get some positive power of $ (1-2 \kappa)^{-1} $ in the right hand side.

\bigskip

\noindent {\it Proof.} The result follows from the form of $ (\overline{P}-z)^{-1}_{\kappa, A,n} $ described in Proposition \ref{algebrecombinatoire} combined with the estimates of Lemma \ref{regulariteSobolevdilatation} and Proposition \ref{fermeturecommutateursprop}. \finpreuve

\bigskip

\begin{coro} \label{corollairecalculfonctionnel} Fix $ n \geq 0 $ integer, $ M \geq 0 $ and  $ \chi \in C_0^{\infty}( \Ra ) $. Then there exists $ C > 0 $ such that for all coefficients $ a_{jk},b_k $ satisfying  1,2 and 3 in Proposition \ref{Propositionsansh}, and for all $ 0 < \kappa \leq 1/4 $,   
the operator $ \chi(\overline{P})_{\kappa, A,n} $ exists  and we have
\begin{eqnarray}
 \big| \big| \chi(\overline{P})_{\kappa, A,n}  \big| \big|_{L^2 \rightarrow L^2} \leq C . \label{sansh}
\end{eqnarray} 
\end{coro}

\noindent {\it Proof.} It is a simple consequence of Proposition \ref{Propositionsansh} and the following Helffer-Sj\"ostrand formula (see for instance \cite{DiSj})
$$ \chi \big( \overline{P} \big) =  \frac{1}{\pi} \int_{\Ca } \bar{\partial} \widetilde{\chi} (z) \big( \overline{P}-z \big)^{-1} L(d z)  , $$
where $ L (d z) $ is the Lebesgue measure on $ \Ca \simeq \Ra^2 $
and $ \widetilde{\chi} \in C_0^{\infty}({\mathbb C}) $ is an almost analytic extension of $ \chi $, ie such that 
$ \overline{\partial} \widetilde{\chi} (z) = {\mathcal O}(|\mbox{Im}(z)|^{\infty}) $ and $ \widetilde{\chi}_{| \Ra} = \chi $. \finpreuve

\bigskip

We shall also need the following proposition.

\begin{prop} There exists $ \epsilon > 0 $ such that for all $ n \geq 0 $ integer and all $ M > 0 $, there exists  $ C > 0 $ such for all coefficients $ a_{jk},b_k $ satisfying 1,2,3 in Proposition \ref{Propositionsansh}, we have
\begin{eqnarray}
 \big| \big| \big(\overline{P}+1 \big)^{-1}_{\kappa, A,n} u \big| \big|_{H^1}  \leq  C ||u||_{H^{-1}} , \label{regularisationresolvente} 
\end{eqnarray}
for all $ u \in L^2 $ and all $ 0 < \kappa \leq 1/4  $. Furthermore, if $ \chi \in C_0^{\infty} $, we have
\begin{eqnarray}
 \big| \big| \chi \big(\overline{P} \big)_{\kappa, A,n} u \big| \big|_{H^1}  \leq  C ||u||_{H^{-1}} . \label{regularisationfonction}
\end{eqnarray}
\end{prop}

\noindent {\it Proof.} It suffices to prove (\ref{regularisationresolvente}) since (\ref{regularisationfonction}) would then follow from (\ref{sansh}) and (\ref{regularisationresolvente}) using the identity
$$  \chi \big(\overline{P} \big)_{\kappa, A,n} =   \big(\overline{P} +1 \big)^{-1}_{\kappa, A,n}  \widetilde{\chi} \big(\overline{P} \big)_{\kappa, A,n}   \big(\overline{P} + 1 \big)_{\kappa, A,n}^{-1} , $$
with $ \widetilde{\chi}(\alpha) = (\alpha+1)^2 \chi (\alpha) $, which is justified by Proposition \ref{Propositionformelle0}. Let us prove (\ref{regularisationresolvente}).
 Using the form of $ \big(\overline{P}+1 \big)^{-1}_{\kappa, A,n} $ given by Proposition \ref{algebrecombinatoire}, the result would follow from the  estimates
\begin{eqnarray}
 \big| \big| \big(\overline{P}+1 \big)^{-1} u \big| \big|_{H^1} & \leq & C ||u||_{H^{-1}}, \qquad u \in L^2 , \label{facile1} \\
\big| \big| \overline{ \mbox{ad}_A^{j} P } v \big| \big|_{H^{-1}} & \leq & C ||v||_{H^1}, \qquad v \in H^2 , \ j \leq n , \label{facile2} \\
\big| \big| (\kappa A+i)^{ - 1} w \big| \big|_{H^{\pm 1}}& \leq & C ||w||_{H^{\pm 1}} , \qquad w \in H^{1 \pm 1} , \label{facile3}
\end{eqnarray}
for some $ C $ independent of the coefficients of $ P $ and $ \kappa $. The estimate (\ref{facile2}) follows from (\ref{1explicitebis}) and (\ref{2explicitebis}). The
estimate (\ref{facile3}) is given in Lemma \ref{regulariteSobolevdilatation} in the $+$ case, and is the adjoint of the $ H^1 \rightarrow H^1 $ bound on $ (\kappa A-i)^{-1} $ in the $-$ case. Finally (\ref{facile1}) follows from the bound
\begin{eqnarray}
 || (\overline{P}+1)^{-1/2} u ||_{H^1} \leq C ||u||_{L^2}, \qquad u \in  L^2 ,  \label{Sobolevimpair1}
\end{eqnarray} 
(and the adjoint one) which follows in a standard fashion from (\ref{faible}). \finpreuve

\section{Elliptic estimates}
\setcounter{equation}{0}

In this section, we prove some elementary elliptic regularity estimates for 
$ \big( \overline{P}+1 \big)_{\kappa, A,n}^{-1} $ (recall Definition \ref{definitionfermeture}). Everywhere we set
$$ r = \bar{r}(d) ,  $$
where $ \bar{r}(d) $ is defined by (\ref{notationregularitemaximale}). We start with the following result.
\begin{prop} \label{continuiteSobolev}  Let $ o \in \{0,1 \} $ and $ s  $ be an integer such that $ 0  \leq s \leq r $.
Then there exists $ C $ such that
\begin{eqnarray}
  || a u ||_{H^{s-o}} \leq C ||a ||_{o,r-o,0} || u ||_{H^s} ,  \label{differentielo}
\end{eqnarray}
for all  $ a \in S^{o,r-o,0} $ and $ u \in {\mathcal S} $. 
\end{prop}

The estimate (\ref{differentielo}) means that the multiplication by $a$ behaves like a differential operator of order $ o $.

\bigskip

\noindent {\it Proof.} We consider first the case when $s=0$ and $ o = 1 $. In this case, the result follows from
$$ || au ||_{H^{-1}} \lesssim || a u ||_{L^{\frac{2d}{d+2}}} \leq || a ||_{L^{d}} || u ||_{L^{2}} , $$
by the H\"older inequality. In the other cases, we have $ s - o \geq 0 $ and we proceed as follows. 
 Observe  that, for any $ 0 \leq k \leq d/2 $,
\begin{eqnarray}
 || \varphi \psi ||_{L^2} \leq || \varphi ||_{L^{\frac{d}{k}}} || \psi ||_{L^{\frac{2d}{d-2k}}} ,   \label{produit}
\end{eqnarray}
by 
the H\"older inequality.  Let $ |\alpha| \leq s-o $. By the Leibniz rule,
$$ || \partial^{\alpha} ( au) ||_{L^2} \leq  C \sum_{\gamma \leq \alpha}  || (\partial^{  \gamma} a ) (\partial^{\alpha-\gamma} u) ||_{L^2} . $$
Since $ u \in H^s $, we have
$$ \partial^{\alpha - \gamma} u \in H^{s - |\alpha| + |\gamma| }  \subset H^{o+|\gamma|}  \subset L^{\frac{2d}{d-2(o+|\gamma|)}}  , $$
the last inclusion being the usual Sobolev embedding (here we use that $ r < d/2 $). By the continuity of this embedding and (\ref{produit}), we have
$$ || (\partial^{ \gamma} a ) (\partial^{\alpha-\gamma} u) ||_{L^2} \leq C || \partial^{\gamma} a ||_{L^{\frac{d}{o+|\gamma|}}} || u ||_{H^{s}} , $$
from which the result follows (recall that $ o+|\gamma| \leq o+|\alpha| \leq s \leq r $). \finpreuve

\bigskip

Using the self-adjointness of $ P $, we obtain the following result for Sobolev spaces of positive or negative order.
\begin{coro} \label{corollairecommutateur} For all integer $n \geq 0 $ and $-r \leq s \leq r$ integer, 
\begin{eqnarray}
|| \emph{ad}_A^n \big( P + \Delta \big) u ||_{H^{s-1}} \lesssim \left( \sum_{jk} || a_{jk} - \delta_{jk} ||_{0,r,n} + ||b_k ||_{1,r-1,n} \right) || u ||_{H^{s+1}},
\end{eqnarray}
for all $ u \in {\mathcal S} $ and all $ a_{jk} \in S^{0,r,n} $, $ b_k \in S^{1,r-1,n} $ such that (\ref{conditionsymetrie}) holds.
\end{coro}

\noindent {\it Proof.} For non negative $s$, the result follows from Proposition \ref{continuiteSobolev} and (\ref{niemecommutateur})-(\ref{formedescoefficients}). For negative $s$, one takes the adjoint since $ i^n \mbox{ad}_A^n (P+\Delta) $ is (formally) self-adjoint. \finpreuve

\bigskip

We next prove the following proposition which will be crucial in Subsection \ref{subsectionsmall}.

\bigskip

\begin{prop} \label{fermetureiso} Fix an integer $ n \geq 0 $ and $ M > 0 $. There exists $ \epsilon > 0 $, $ \kappa_0 > 0 $ and $ C > 0 $ such that
if 
\begin{enumerate}
\item{ $ ||a_{jk} - \delta_{jk} ||_{0,r,0} + ||b_k ||_{1,r-1,0} < \epsilon  $, }
\item{ $  ||a_{jk} ||_{0,r,n} + || b_k ||_{1,r,n} \leq M$,}
\item{ $ 0 < \kappa \leq \kappa_0 $,}
\item{ $ - r \leq s \leq r $ integer,}
\end{enumerate}
then the operator
$$ (\kappa A+i)^{-n}(P+1) (\kappa A+i)^n : {\mathcal S} \rightarrow H^{s-1} \cap L^2  ,$$
has a bounded closure  $ H^{s + 1} \rightarrow H^{s-1} $ denoted by $ (P+1)_{\kappa, A,n,s} $ which is an isomorphism between $ H^{s+1} $ and $ H^{s-1} $
and such that
\begin{eqnarray}
 ||u||_{H^{s+1}} / C \leq || (P+1)_{\kappa, A,n,s} u ||_{H^{s-1}} \leq C || u ||_{ H^{s +1} } , \qquad u \in H^{s+1} . \label{avecladensite} 
\end{eqnarray}
Furthermore, with the notation of Definition \ref{definitionfermeture},
\begin{eqnarray}
 \big| \big| (\overline{P}+1)^{-1}_{\kappa, A,n}  u \big| \big|_{H^{s+1}} \leq C || u ||_{H^{s-1}}, \qquad u \in L^2 \cap H^{s-1} . \label{bonnejustificationformelle}
\end{eqnarray}
\end{prop}

\noindent {\it Proof of Proposition.} Observe first  that
$$ (\kappa A+i)^{-n} P (\kappa A+i)^n - P  $$
is a linear combination, with coefficients which are universal constants, of operators of the form
\begin{eqnarray*}
 \kappa^m (\kappa A+i)^{-m_1} \mbox{ad}_A^m (P)  , \qquad m_1 \geq 0, \ 1 \leq m \leq n .
\end{eqnarray*}
Therefore, using Corollary \ref{corollairecommutateur} and Lemma \ref{regulariteSobolevdilatation} (with $ \kappa(|r|+1)<1/2 $) we have
$$ || (\kappa A+i)^{-n} P (\kappa A+i)^n u + \Delta u  ||_{H^{s-1}} \lesssim  \left(\epsilon + \kappa M  \right) ||u||_{H^{s+1}}  . $$
By choosing $ \epsilon $ and $ \kappa $ small enough, we obtain the existence of the closure $ (P+1)_{\kappa, A,n,s} $ and the fact that it is close to $ 1 - \Delta $ in the $ H^{s+1} \rightarrow H^{s-1} $ topology, hence is an isomorphism. We also get (\ref{avecladensite}).
 To prove (\ref{bonnejustificationformelle}) it suffices to show that
 \begin{eqnarray}
  (\overline{P}+1)^{-1}_{\kappa, A,n} =   \left( (P+1)_{\kappa, A,n,s} \right)^{-1} \qquad \mbox{on} \ \  L^2 \cap H^{s-1} , \label{yielding}
 \end{eqnarray} 
and then use the lower bound in (\ref{avecladensite}). One sees that (\ref{yielding}) holds by checking that
\begin{eqnarray}
   (\overline{P}+1)^{-1}_{\kappa, A,n}  (P+1)_{\kappa, A,n,s} w = w , \qquad s - 1 \geq 0 ,  \nonumber \\
     (P+1)_{\kappa, A,n,s} (\overline{P}+1)^{-1}_{\kappa, A,n} w = w , \qquad s-1 < 0 , \nonumber
\end{eqnarray}
for all $w \in {\mathcal S} $. This follows from Lemma \ref{Lemme1} by approaching $ u = (P+1) (\kappa A+i)^n w $ in $ H^{s-1} $ (hence in $ L^2 $) in the first case, and  $ u = (\overline{P}+1)^{-1}(\kappa A+i)^n w$  in $ H^2 $ (hence in $ H^{s+1} $) in the second case.
\finpreuve

\section{Resolvent estimates}
\setcounter{equation}{0}
The purpose of this section is to prove Theorem \ref{estimationprincipale}. The latter will be divided into two steps. In Subsection \ref{subsectionsmall}, 
we shall prove  resolvent estimates for operators of the form (\ref{formedePdetail}) which are small perturbations of $ - \Delta $, using a scale invariant analysis.
In Subsection \ref{nonsmallsubsection}, we will prove Theorem  \ref{estimationprincipale} by combining a compactness argument and the estimates of Subsection \ref{subsectionsmall}, by reducing $ - \Delta_G $ to a compactly supported perturbation of an operator of the form (\ref{formedePdetail}).

\subsection{Small perturbations} \label{subsectionsmall}
Throughout this subsection, $ P $ denotes an operator of the form (\ref{formedePdetail}) and, as before, $ \overline{P} $ denotes its $ H^2 \rightarrow L^2 $ closure which is selfadjoint on $ L^2 $ with domain $ H^2 $. We shall basically prove weighted estimates on $ (\overline{P}-z)^{-n} $ seen as an operator from $ H^{-n} \rightarrow H^n $. The first step is to get $ L^2 \rightarrow L^2 $ estimates and is the purpose of the following proposition.

\begin{prop}[Jensen-Mourre-Perry estimates \cite{JeMoPe}] \label{propJeMoPe} There exists $ \epsilon > 0 $ such that, for all integer $ N \geq 0 $, all $ M \geq 0 $ and all relatively compact interval $ I \Subset (0,\infty) $, there exists $ C > 0 $ such that, for all $ n \leq N $,
$$ \big| \big| (A+i)^{-n} \big( \overline{P} - z \big)^{-n} (A+i)^{-n} \big| \big|_{L^2 \rightarrow L^2} \leq C, \qquad \emph{Re}(z) \in I, \ \emph{Im}(z) \ne 0 , $$
for all coefficients $ a_{jk} \in S^{0,1,N+1}$, $ b_k \in S^{1,0,N+1} $ such that
 \begin{enumerate}
 \item{(\ref{conditionsymetrie}) is satisfied}
 \item{$ \sum_{j,k} ||a_{jk} - \delta_{jk}||_{0,1,1} + ||b_k||_{1,0,1} < \epsilon $ }
 \item{$ \sum_{j,k} ||a_{jk} ||_{0,1,N+1} + ||b_k||_{1,0,N+1} \leq M $.}
\end{enumerate}
\end{prop}

The latter result follows by tracking the uniform dependence of the estimates with respect the coefficients of $ P $ in the proofs of \cite{JeMoPe}. We simply point out that the smallness of $ ||a_{jk} - \delta_{jk} ||_{0,1,1} $ and $ ||b_k ||_{1,0,1} $ guarantees the positive commutator estimate
$$ \chi (\overline{P}) i \overline{[P,A]} \chi (\overline{P}) \geq \frac{\inf {\rm supp}(\chi)}{2} \chi^2 (\overline{P}) , $$
for any $ \chi \in C_0^{\infty}(\Ra^+) $, which follows from (\ref{commutateurexplicitenouveau}). The other ingredient is the uniform $ \overline{P} $ boundedness estimate
(\ref{Pborneuniforme}).
 
\bigskip

In the sequel, we shall use  the notation
\begin{eqnarray}
 a_{\tau} (x) := a (e^{\tau} x), \qquad \tau \in \Ra, \ x \in \Ra^d ,  \label{actiondilatation}
\end{eqnarray}
namely
$$ a_{\tau} = e^{i \tau A} a e^{-i\tau A} . $$
Here is the main property of the spaces $ S^{o,r-o,N} $.
\begin{prop}[Scaling homogeneity] \label{scaleinvariance} Let $ o \in \{0,1 \} $ and $ o \leq r \leq \bar{r}(d) $. For all $ a \in S^{o,r-o,N} $ and $ \tau \in \Ra $,
$$ e^{o \tau } ||a_{\tau} ||_{o,r-o,N} = || a ||_{o,r-o,N} . $$
\end{prop}

\noindent {\it Proof.} Observe first that
\begin{eqnarray}   \left( (x \cdot \nabla )^n a \right)_{\tau} =  (x \cdot \nabla )^n (a_{\tau} )   , \label{invariancedilatation}
\end{eqnarray}
either by a trivial direct computation, or by remarking that dilations commute with their generator.
Then 
$$  (x \cdot \nabla)^n ( \partial^{\beta} a_{\tau})   = e^{\tau |\beta|} \left(  (x \cdot \nabla)^n \partial^{\beta} a \right)_{\tau} , $$
and we see that
$$ e^{\tau o} || (x \cdot \nabla)^n\partial^{\beta} (a_\tau) ||_{L^{\frac{d}{|\beta|+o}}}  =  || (x \cdot \nabla)^n\partial^{\beta} a ||_{L^{\frac{d}{|\beta|+o}}}, $$
by an elementary change of variable in the integral when $ |\beta| + o \ne 0 $, and trivially if $ |\beta| + o = 0 $.  \finpreuve

\bigskip
 
We are now ready to prove the following theorem which is our main technical result. 

\begin{theo} \label{theoremeprincipalpetit} Let $ N := \bar{r}(d)+1 $. Fix a constant $ M > 0 $. Then, there exist $ \epsilon > 0 $ and $ \kappa > 0$ such that, for all $ a_{jk} \in S^{0,\bar{r}(d),N+1} $,
$ b_k \in S^{1,\bar{r}(d)-1,N+1} $ such that
\begin{enumerate}
\item{ (\ref{conditionsymetrie}) is satisfied}
\item{ $ \sum_{j,k} ||a_{jk} - \delta_{jk} ||_{0,1,1} + ||b_k||_{1,0,1} < \epsilon $,}
\item{ $ \sum_{j,k} ||a_{jk} - \delta_{jk} ||_{0,\bar{r}(d),0} + ||b_k||_{1,\bar{r}(d)-1,0} < \epsilon $,}
\item{ $ \sum_{j,k} ||a_{jk}  ||_{0,\bar{r}(d),N+1} + ||b_k||_{1,\bar{r}(d)-1,N+1} < M $, }
\end{enumerate}
we have the following estimates,
\begin{itemize}
\item{ if $ 1 \leq n \leq \bar{r}(d) $
\begin{eqnarray}
\big| \big| (\kappa A+i)^{-n} ( \overline{P} - z)^{-n} (\kappa A-i)^{-n} \varphi \big| \big|_{L^{q(n)}} \leq C  ||\varphi ||_{L^{p(n)}}, \label{casgeneral}
\end{eqnarray}}
\item{if $ n = N $ and $ d $ is odd,
\begin{eqnarray}
\big| \big| (\kappa A+i)^{-N} ( \overline{P}-  z)^{-N} (\kappa A-i)^{-N} \varphi \big| \big|_{L^{\infty}} \leq C \emph{Re}(z)^{-1/2} ||\varphi ||_{L^{1}}, \label{casimpair}
\end{eqnarray}}
\item{if $ n= N $ and $d$ is  even, then for all $  2d < q < \infty $,  
\begin{eqnarray}
\big| \big| (\kappa A+i)^{-N} (\overline{P}-  z)^{-N} (\kappa A-i)^{-N} \varphi \big| \big|_{L^{q}} \leq C_q  \emph{Re}(z)^{- \frac{2d}{q}} ||\varphi ||_{L^{\frac{q}{q-1}}} ,
\label{caspair}
\end{eqnarray}}
\end{itemize}
all these estimates holding for
$$ \varphi \in {\mathcal S}, \qquad \emph{Re}(z)> 0, \qquad \emph{Im}(z) \ne 0 . $$
\end{theo}

This result can be viewed as a version of Theorem \ref{estimationprincipaleLp} for small perturbations of the flat metric. Notice that the coefficients of the perturbation are taken in the classes $ S^{o,r-o,N} $, which is a more general condition than being in $ S^{-\rho-o} $, as we shall see in Proposition \ref{classe}. Note also that Theorem \ref{theoremeprincipalpetit} holds for $ \mbox{Re}(z) $ small, which is its main interest, but actually for all $ \mbox{Re}(z) > 0 $ hence for large ones too. 

\bigskip

\noindent {\it Proof.} It is based on an scaling argument. Let $ \lambda = \mbox{Re}(z) $ and write
$$ P - z = \lambda \big( \lambda^{-1} P - 1 - i \delta \big) , $$
where $ \lambda \delta = \mbox{Im}(z) $. Then, by setting $$ \lambda^{-1/2} = e^{\tau} , $$
and
$$  P_{\tau} = a_{jk,\tau} (x) D_j D_k + e^{\tau }b_{k,\tau}(x)D_k , $$
where we use the notation (\ref{actiondilatation}), we have
$$ \lambda^{-1} P = e^{-i \tau A} P_{\tau} e^{i \tau A } , $$
and thus
\begin{eqnarray}
 \big(\overline{P}-z \big)^{-1} = \lambda^{-1} e^{- i \tau A} \big( \overline{P}_{\tau} - 1 - i \delta \big)^{-1} e^{i \tau A} . \label{autiliseravantscaling}
\end{eqnarray}
By Proposition \ref{scaleinvariance}, the conditions {\it 2}, {\it 3} and {\it 4} hold for $ a_{jk,\tau} $ and $ e^{\tau} b_{k,\tau} $, uniformly with respect to $ \tau \in \Ra $.
In particular, using Proposition \ref{propJeMoPe}, we have
\begin{eqnarray}
 \big| \big| (A+i)^{-n} \big( \overline{P}_{\tau} - 1 - i \delta \big)^{-n} (A-i)^{-n} \big| \big|_{L^2 \rightarrow L^2} \leq C_M , \label{resolventeL2L2}
\end{eqnarray}
 for all $ \tau \in \Ra$, $ \delta \in \Ra \setminus 0 $ and $ \ 1 \leq n \leq N $. Since this is an $ L^2 \rightarrow L^2 $ estimate, $ (A\pm i)^{-n} $ can be replaced by $ (\kappa A \pm i)^{-n} $ therein, for any $ \kappa > 0$, up to the replacement of $C_M$ by a $ \kappa $ dependent constant. This will be useful to consider $ H^{-n} \rightarrow H^{n} $ estimates as follows. Introduce $ \chi \in C_0^{\infty} (\Ra) $ which is real valued and equal to $ 1 $ near $1$. We then split the resolvent as
\begin{eqnarray}
 \big( \overline{P}_{\tau} - 1 - i \delta \big)^{-n} & = &\big( \overline{P}_{\tau} - 1 - i \delta \big)^{-n} (1 - \chi^2) \big( \overline{P}_{\tau} \big)
+ \chi \big( \overline{P}_{\tau} \big) \big( \overline{P}_{\tau} - 1 - i \delta \big)^{-n} \chi \big( \overline{P}_{\tau} \big) , \nonumber \\
 & = & \mbox{I}(\tau,\delta) + \mbox{II}(\tau,\delta) .
\label{decompositionspectralesimple}
\end{eqnarray}
We consider first $ \mbox{I}(\tau,\delta) $. By setting  $ \Phi_{\delta} (\alpha) = (1-\chi^2 (\alpha)) (\alpha + 1)^n / (\alpha - 1 - i \delta)^n $, we  can write
$$ \mbox{I}(\tau,\delta) = \big( \overline{P}_{\tau} + 1  \big)^{-n/2} \Phi_{\alpha} (\overline{P}_{\tau}) \big( \overline{P}_{\tau} + 1   \big)^{-n/2} . $$
Since $ \Phi_{\delta} $ is bounded in $ L^{\infty} ([0,\infty)_{\alpha}) $ as $ \delta  $ varies, the Spectral Theorem yields
$$ ||\Phi_{\delta} (\overline{P}_{\tau}) ||_{L^2 \rightarrow L^2} \leq C , \qquad\tau \in \Ra, \ \ \delta \in \Ra \setminus 0 . $$
On the other hand, using (\ref{bonnejustificationformelle}), and also (\ref{Sobolevimpair1}) if $n$ is odd, we have
$$ || ( \overline{P}_{\tau} + 1 )^{-n/2} \psi ||_{H^n} \leq C || \psi ||_{L^2}, \qquad \psi \in L^2, \ \tau \in \Ra . $$
We also have the dual $ H^{-n} \rightarrow L^2 $ bound and we conclude, using Lemma \ref{regulariteSobolevdilatation},  that
$$ \big| \big| (\kappa A+i)^{-n} \mbox{I}(\tau,\delta) (\kappa A-i)^{-n} \psi \big| \big|_{H^n} 
\leq C || \psi ||_{H^{-n}}, $$
for all $ \psi \in L^2 $, $ \tau \in \Ra $ and $ \delta \ne 0 $. We next consider the second term of (\ref{decompositionspectralesimple}). By Proposition \ref{Propositionformelle0} and Corollary \ref{corollairecalculfonctionnel}, we can write
$$ (\kappa A+i)^{-n} \mbox{II}(\tau,\delta) (\kappa A-i)^{-n} = \chi \big( \overline{P}_{\tau} \big)_{\kappa, A,n} (\kappa A+i)^{-n} \big( \overline{P}_{\tau} - 1 - i \delta \big)^{-n} (\kappa A-i)^{-n}  \chi \big( \overline{P}_{\tau} \big)_{\kappa, A,n}^* . $$
We then observe that we have the estimate
\begin{eqnarray}
 || \chi \big( \overline{P}_{\tau} \big)_{\kappa, A,n} ||_{L^2 \rightarrow H^n} \leq C_n, \qquad \tau \in \Ra . \label{rereference}
\end{eqnarray}
The latter is obtained by writting $ \chi (\alpha) = (\alpha + 1)^{-[n/2]} \psi (\alpha)  $ with $ [n/2] $ the integer part of $ n/2 $ so that
$$  \chi \big( \overline{P}_{\tau} \big)_{\kappa, A,n} =   \big( \overline{P}_{\tau} + 1 \big)_{\kappa, A,n}^{-[n/2]}  \psi \big( \overline{P}_{\tau} \big)_{\kappa, A,n} : L^2 \rightarrow H^n ,$$ 
 by   (\ref{bonnejustificationformelle}) and Proposition \ref{Propositionformelle0} if $n$ is even or (\ref{bonnejustificationformelle}) and (\ref{regularisationfonction}) if $n$ is odd. Similarly, we have a $ H^{-n} \rightarrow L^2 $ bound for $ \chi \big( \overline{P}_{\tau} \big)_{\kappa, A,n}^*$.  Thus, using the $ L^2 \rightarrow L^2 $ bound (\ref{resolventeL2L2}),
we deduce that
$$ \big| \big| (\kappa A+i)^{-n} \mbox{II}(\tau,\delta) (\kappa A-i)^{-n} \psi \big| \big|_{H^n} 
\leq C || \psi ||_{H^{-n}}, $$
and conclude that
$$  \big| \big| (\kappa A+i)^{-n} \big( \overline{P}_{\tau} - 1 - i \delta \big)^{-n} (\kappa A-i)^{-n} \psi \big| \big|_{H^n} 
\leq C || \psi ||_{H^{-n}}, $$
for all $ \psi \in L^2 $, $ \tau \in \Ra $ and $ \delta \ne 0 $. In terms of $ \overline{P} $ the latter reads
\begin{eqnarray}
  \big| \big| e^{i \tau A }(\kappa A+i)^{-n} \big( \overline{P} - z \delta \big)^{-n} (\kappa A-i)^{-n}  \varphi \big| \big|_{H^n} , \label{redaction}
\leq C \lambda^{-n} || e^{i \tau A } \varphi ||_{H^{-n}}, 
\end{eqnarray}
where we recall that $ \lambda = \mbox{Re}(z) $. We can get rid of the negative powers of $ \lambda $ as follows.
If $ n \leq N-1 $, we have on one hand the Sobolev embeddings 
$$ L^{p(n)} \subset H^{-n}, \qquad H^n \subset L^{q(n)} . $$
On the other hand, using (\ref{pourlescalingfinal}) and the fact that $ \lambda^{-1} = e^{2 \tau} $, we have  
\begin{eqnarray}
 \lambda^{-n} \ || e^{i \tau A} ||_{L^{p(n)} \rightarrow L^{p(n)}} || e^{i \tau A} ||_{L^{q(n)} \rightarrow L^{q(n)}}^{-1} = 1 . \label{calculscaling}
\end{eqnarray}
Thus, by turning (\ref{redaction}) into a $ L^{ p(n)} \rightarrow L^{q(n)} $ estimate and by using (\ref{calculscaling}), we obtain (\ref{casgeneral}). 
If $ n  = N $, the same argument applies using the Sobolev embeddings with
$$ q(N)= \infty \ \ \mbox{if} \ d \mbox{\ is odd}, \qquad q(N) = q \ \mbox{with an arbitrary} \ q > 2d \ \ \mbox{if} \ d \ \mbox{is even} , $$
the only difference being that the left hand side of (\ref{calculscaling}) becomes either $ \lambda^{-1/2} $ or $ \lambda^{-2d/q} $. \finpreuve

\bigskip

To apply Theorem \ref{theoremeprincipalpetit} to perturbations of the Laplacian with coefficients in $ S^{-\rho} $, we need the following result.

\begin{prop}[Symbol classes embeddings] \label{classe} For all $ 1 \leq r \leq \bar{r}(d)  $, $ N \geq 0 $ integers and  $ \mu > 0 $ real, we have the continuous embeddings
\begin{eqnarray}
S^{-\mu-1} & \subset & S^{1,r-1,N}, \nonumber \\
 S^{-\mu} & \subset & S^{0,r,N} . \nonumber
\end{eqnarray}
\end{prop}

These embeddings are very convenient since the seminorms of the spaces $ S^{-\mu} $ behave badly under scaling, unlike $ S^{o,r-o,N} $ by Proposition \ref{scaleinvariance}.

\bigskip

\noindent {\it Proof.} We note first that
 $$ \scal{x}^{-\mu-|\beta|-o} \in L^{\frac{d}{|\beta|+o}}  . $$
Furthermore, by an elementary induction, one checks that  $ (x \cdot \nabla)^n $ is a linear combination of $ x^{\alpha} \partial^{\alpha} $ with $ | \alpha | \leq n $. Therefore,  we have the estimates
$$ || (x \cdot \nabla)^n \partial^{\beta} a ||_{L^{\frac{d}{|\beta|+ o}}}  \leq C || \scal{x}^{\mu + o + |\beta|} (x \cdot \nabla)^n \partial^{\beta} a ||_{L^\infty} \leq C \max_{|\alpha| \leq n} || \scal{x}^{\mu + o + |\beta|+|\alpha|} \partial^{\alpha+\beta} a ||_{L^\infty}  , $$
which lead easily to the result. \finpreuve

\bigskip

By this proposition, we see that Theorem \ref{theoremeprincipalpetit} holds if the coefficients of $ P $ are such that $ a_{jk}-\delta_{jk} $ and $ b_k $ are small enough respectively in $ S^{-\rho} $ and $ S^{-1-\rho} $. We may also replace the weights $ (\kappa A \pm i)^{-1} $ by powers of $ \scal{x}^{-1} $ according to a classical procedure. This is the purpose of the following.
\begin{coro} \label{corollairepratique1} Assume (\ref{conditionsymetrie}) and that $ a_{jk} - \delta_{jk} \in S^{-\rho} $ and  $ b_k \in S^{-1-\rho} $. Assume also that
\begin{eqnarray}
 \big| \scal{x}^{\rho + |\alpha|} \partial^{\alpha} (a_{jk}(x) - \delta_{jk}) \big| + \big| \scal{x}^{1+\rho+|\alpha|} \partial^{\alpha}b_k (x) \big| \leq \epsilon ,
 \label{conditiondepetitesse} 
\end{eqnarray}
for $ |\alpha| \leq \bar{r}(d) + 1 $. If $ \epsilon $ is small enough, then for $ 1 \leq n \leq N := \bar{r}(d) + 1 $
\begin{itemize}
\item{ if $ 1 \leq n \leq \bar{r}(d) $
\begin{eqnarray}
\big| \big| \scal{x}^{-n} ( \overline{P} - z)^{-n} \scal{x}^{-n} \varphi \big| \big|_{L^{q(n)}} \leq C  ||\varphi ||_{L^{p(n)}}, \nonumber
\end{eqnarray}}
\item{if $ n = N $ and $ d $ is odd,
\begin{eqnarray}
\big| \big| \scal{x}^{-N} ( \overline{P}-  z)^{-N} \scal{x}^{-N} \varphi \big| \big|_{L^{\infty}} \leq C \emph{Re}(z)^{-1/2} ||\varphi ||_{L^{1}}, \nonumber
\end{eqnarray}}
\item{if $ n= N $ and $d$ is  even, then for all $  2d < q < \infty $,  
\begin{eqnarray}
\big| \big| \scal{x}^{-N} (\overline{P}-  z)^{-N} \scal{x}^{-N} \varphi \big| \big|_{L^{q}} \leq C_q  \emph{Re}(z)^{- \frac{2d}{q}} ||\varphi ||_{L^{\frac{q}{q-1}}} , \nonumber
\end{eqnarray}}
\end{itemize}
 for
$$ \varphi \in {\mathcal S}, \qquad 0 < \emph{Re}(z)< 1, \qquad \emph{Im}(z) \ne 0 . $$
In particular, we may replace all $ L^{p} , L^{q} $ spaces above by $ L^2 $ if we change $ \scal{x}^{-n} $ into $ \scal{x}^{-2n - \varepsilon} $,
  for any $ \varepsilon > 0 $.
\end{coro} 
 
\medskip 
 
\noindent {\it Proof.} This kind of result is standard so we briefly recall the proof. Note first that (\ref{conditiondepetitesse}) implies that items 2 and 3 of Theorem \ref{theoremeprincipalpetit} are satisfied. Fix $ \chi \in C_0^{\infty} (\Ra) $ such that $ \chi \equiv 1 $ near $ [0,1] $. By the Spectral Theorem and elliptic regularity
$$ (1 - \chi (\overline{P})) (\overline{P}-z)^{-n} , $$
maps $ H^{-n} $ to $ H^n $, uniformly with respect to $z$, hence the appropriate Lebesgue spaces to their duals by Sobolev embeddings. Thus, it suffices to consider
$$ \scal{x}^{ - n } \big( \overline{P} -z \big)^{-n} \chi (\overline{P}) \scal{x}^{ -  n} , \qquad \mbox{Re}(z) \in (0,1), \ \mbox{Im}(z) \ne 0 . $$
By possibly choosing $ \chi $ of the form $ \varphi^2 $ the result follows by writting
 $$ \scal{x}^{ - n } \varphi (\overline{P}) =  \left( \scal{x}^{-n} \varphi (\overline{P}) (\kappa A+i)^n \right) (\kappa A+i)^{-n} ; $$
where we observe that, for all $ q \in [1,\infty] $,
$$ \scal{x}^{-n} \varphi (\overline{P}) (\kappa A+i)^n : L^{q} \rightarrow L^{q} , $$
for it is a pseudo-differential operator with symbol in $ S^{-\infty} $ (see {\it e.g.} \cite{BTJFA}). The replacement of $ L^q $ spaces by $ L^2 $ after the replacement of $ \scal{x}^{-n} $ by $ \scal{x}^{-2n-\epsilon} $ follows from
\begin{eqnarray}
|| \scal{x}^{-n - \varepsilon}  v ||_{L^2} \leq C || v ||_{ L^{q(n) } }  , \label{reutiliseunefois}
\end{eqnarray}
by the H\"older inequality (note that this works even for $n = N$ and $ q(N) \in ( 2d , \infty]  $). \finpreuve

\bigskip

In the next paragraph, we will also need the following result for $ \mbox{Re}(z)<0 $.
\begin{prop} \label{estimeecotenegatif} Under the same assumptions as in Theorem \ref{theoremeprincipalpetit},
we have the following estimates:
\begin{itemize}
\item{ if $ 1 \leq n \leq \bar{r}(d) $
\begin{eqnarray}
\big| \big|  ( \overline{P} - z)^{-n}  \varphi \big| \big|_{L^{q(n)}} \leq C  ||\varphi ||_{L^{p(n)}}, \label{casgeneral2}
\end{eqnarray}}
\item{if $ N = \bar{r}(d)+1 $ and $ d $ is odd,
\begin{eqnarray}
\big| \big|  ( \overline{P}-  z)^{-N}  \varphi \big| \big|_{L^{\infty}} \leq C |\emph{Re}(z)|^{-1/2} ||\varphi ||_{L^{1}},  \label{casimpair2}
\end{eqnarray}}
\item{if $ N = \bar{r}(d)+1 $ and $d$ is  even, then for all $  2d < q < \infty $,  
\begin{eqnarray}
\big| \big|  (\overline{P}-  z)^{-N}  \varphi \big| \big|_{L^{q}} \leq C_q |\emph{Re}(z)|^{- \frac{2d}{q}} ||\varphi ||_{L^{\frac{q}{q-1}}} ,
\label{caspair2}
\end{eqnarray}}
\end{itemize}
all these estimates holding for
$$ \varphi \in {\mathcal S}, \qquad \emph{Re}(z) < 0, \qquad \emph{Im}(z) \ne 0 . $$
\end{prop}

\bigskip

\noindent {\bf Remark} Since $ (\kappa A \pm i)^{-1} $ preserve all $ L^p $ spaces for $\kappa$ small enough, we may replace $ (\overline{P}-z)^{-n} $ by 
$$ (\kappa A+i)^{-n} (\overline{P} - z)^{-n} (\kappa A-i)^{-n} , $$ for $ 1 \leq n \leq N $, in the estimates (\ref{casgeneral2}), (\ref{casimpair2}) and (\ref{caspair2}). In particular, this shows that the estimates of Theorem \ref{theoremeprincipalpetit} 
actually hold for $ \mbox{Re}(z) \in \Ra $. Also, using (\ref{reutiliseunefois}), we may clearly turn all the estimates of Proposition \ref{estimeecotenegatif} into $ L^2 \rightarrow L^2 $ estimates with weights.

\bigskip

\noindent {\it Proof of Proposition \ref{estimeecotenegatif}.} It is based on the same scaling argument as  the proof of Theorem \ref{theoremeprincipalpetit}, from which we borrow the notation. We write $ \mbox{Re}(z) = - \lambda $ so that
$$ \big(\overline{P}-z \big)^{-1} =  \lambda^{-1} \big( \lambda^{-1} P + 1 + i \delta \big)^{-1} $$
the left hand side of which we write as
$$ \lambda^{-1} e^{i \tau A} \big( P_{\tau} + 1 + i \delta \big)^{-1} e^{-i \tau A} . $$
We may then write
$$ \big( P_{\tau} + 1  \big)^{-n/2} \left( \frac{ P_{\tau} + 1  }{ P_{\tau} + 1 + i \delta } \right)^n \big( P_{\tau} + 1  \big)^{-n/2} . $$
The latter is bounded from $ H^{-n} $ to $ H^n $, uniformly with respect to $ \tau $ and $ \delta $ (see (\ref{rereference})) and we conclude as in the proof of Theorem \ref{theoremeprincipalpetit}. \finpreuve

\subsection{Non small perturbations} \label{nonsmallsubsection}

The purpose of this paragraph is to prove Theorems \ref{estimationprincipaleLp} and \ref{estimationprincipale}. We shall actually prove Theorem \ref{estimationprincipale} first and then Theorem \ref{estimationprincipaleLp}.

We start by doing some reductions.

We first choose suitable coordinates on $ \Ra^d $ such that we may assume that $ \mbox{det} \ G (x) = 1 $ outside a compact set. This is explained in Appendix \ref{coordonneesadaptees}. We next conjugate in the usual way our Laplacian to get an   operator which is self-adjoint with respect to the Lebesgue measure:
the map $ u \mapsto \mbox{det} \ G (x)^{1/4} u $ is  unitary from $ L^2 (\Ra^d , d_G x) $ onto $ L^2 (\Ra^d, dx) $ so $ - \Delta_G $ is unitarily equivalent to the operator
\begin{eqnarray}
 P = - \mbox{det}\ G(x)^{-1/4} \frac{\partial}{\partial x_j} \left( \mbox{det}\ G(x)^{1/2}  G_{jk}(x)  \frac{\partial}{\partial x_k} \right) \mbox{det} \ G (x)^{-1/4} , 
 \nonumber
\end{eqnarray}
which has a self-adjoint closure $ \overline{P} $, with domain $ H^2 $.

One may then clearly write
$$ P = P_0 + W , $$
with
$$ P_0 = a_{jk} (x) D_j D_k + b_k (x) D_k , $$
and
\begin{eqnarray}
 W = \chi_{jk}(x)D_j D_k  + \theta_k (x) D_k + V (x) . \label{decompositioncompacte}
\end{eqnarray}
such that
\begin{eqnarray}
 a_{jk} - \delta_{jk} \ \mbox{is small enough in} \ S^{-\frac{\rho}{2}}, \qquad b_k \ \mbox{is small enough in} \ S^{-1-\frac{\rho}{2}} , \label{precision} 
\end{eqnarray}
and
$$ \chi_{jk}, \theta_k , V \in C_0^{\infty} . $$
By small enough, we mean in (\ref{precision}) that we may assume that the estimates of Corollary \ref{corollairepratique1} and Proposition \ref{estimeecotenegatif} hold for $ \overline{P}_0 $.
Here and in the sequel we denote by $ \overline{P}_0 $, $ \overline{P} $ and $ \overline{W} $ the $ H^2 \rightarrow L^2 $ closures of the corresponding differential operators which are a priori defined on $ {\mathcal S} $. In particular, $ \overline{P}_0 $ and $ \overline{P} $ are self-adjoint with domain $ H^2 $ and, by unitary equivalence with $ - \Delta_G $, we have
$$ \overline{P} \geq 0 \qquad \mbox{and} \qquad  0 \ \mbox{is not an eigenvalue of} \ \overline{P} . $$

\bigskip
By the Spectral Theorem, it is sufficient to prove Theorem \ref{estimationprincipale} with $ (\overline{P}-z)^{-n} $ replaced by
\begin{eqnarray}
 R_{\psi}^n (z) := \psi \big( \overline{P} \big) \big( \overline{P} - z \big)^{-n} , \label{Resolvantepsi}
\end{eqnarray}
for some $$ \psi \in C_0^{\infty} (\Ra), \qquad \psi \equiv 1 \ \ \mbox{near} \ 0 . $$
It is also convenient to introduce  $ \Psi \in C_0^{\infty} (\Ra) $ such that 
$$ \Psi \psi = \psi . $$
Both are chosen with values in $ [0,1] $. 
\begin{prop} Let us set
 $$ S_{\Psi}(z) =  \overline{W} \big( \overline{P}_0-z \big)^{-1}  \Psi \big( \overline{P} \big) , $$
 and
 $$ B^1 (z) = \Psi \big(\overline{P} \big) \big(\overline{P}_0-z \big)^{-1} \psi \big(\overline{P} \big) - \Psi \big(\overline{P} \big) \big(\overline{P}_0-z \big)^{-1} \overline{W} \psi \big(
 \overline{P} \big) \big(\overline{P}_0-z \big)^{-1} \Psi \big(\overline{P}\big) . $$
Then
\begin{eqnarray}
  R_{\psi}^1(z) = B^1 (z) + S_{\Psi}(\bar{z})^* R_{\psi}^1(z) S_{\Psi}(z) , \label{pourNeumann}
\end{eqnarray}
for all $ z \in \Ca \setminus \Ra $.
\end{prop}

\noindent {\it Proof.} It is based on the resolvent identity, namely
\begin{eqnarray}
\big(\overline{P} - z \big)^{-1} & = & \big(\overline{P}_0-z \big)^{-1} - \big(\overline{P}_0 - z \big)^{-1} \overline{W} \big(\overline{P}-z \big)^{-1}, \label{identitedelaresolvente} \\
 & = & \big( \overline{P}_0-z \big)^{-1} -  \big(\overline{P} - z \big)^{-1} \overline{W} \big(\overline{P}_0-z \big)^{-1}. \label{identitedelaresolvente2}
\end{eqnarray}
The identity (\ref{pourNeumann}) is obtained by applying first $ \psi (\overline{P}) $ to the right of both sides of (\ref{identitedelaresolvente}), then by inserting (\ref{identitedelaresolvente2})
on the right hand side of the resulting identity and finally by applying $ \Psi (\overline{P})$ to the left and right. \finpreuve

\bigskip

Our strategy is to   show that one can make $ S_{\Psi}(z) $ small enough (in operator norm on suitable weighted $ L^2 $ spaces) by choosing $ \Psi $ (and hence $ \psi $) with a small enough support around $ 0 $ and by choosing $z$ close enough to $ 0 $. To this end, we denote
\begin{eqnarray}
 z = \lambda + i \epsilon , \label{notationz}
\end{eqnarray} 
and introduce the decomposition
\begin{eqnarray}
 S_{\Psi}(z) & = &  \overline{W} \big( \overline{P}_0- i \epsilon \big)^{-1}  \Psi \big( \overline{P} \big) \ + \
 \overline{W} \left( \big( \overline{P}_0-z \big)^{-1} - 
\big( \overline{P}_0- i \epsilon \big)^{-1} \right)  \Psi \big( \overline{P} \big) . \label{Neumann2}
\end{eqnarray}
\begin{prop} \label{absorbe1} Fix $ M > 0 $ and $ \nu > 1 $. If $ \emph{supp}(\Psi) $ is contained in a sufficiently small neighborhood of $ 0 $, then
$$ \big| \big| \scal{x}^M \overline{W} \big( \overline{P}_0- i \epsilon \big)^{-1}  \Psi \big( \overline{P} \big) \scal{x}^{-\nu} \big| \big|_{L^2 \rightarrow L^2} \leq \frac{1}{4}, \qquad \epsilon > 0 . $$
\end{prop}

\noindent {\it Proof.} It suffices to show that, for some $ \delta > 0 $ as small as we want,
\begin{eqnarray}
 \big| \big| \scal{x}^M \overline{W} \big( \overline{P}_0- i \epsilon \big)^{-1}  \Psi \big( \overline{P} \big) \scal{x}^{-\nu} \big| \big|_{L^2 \rightarrow L^2} \leq C \big| \big| \Psi (\overline{P}) \scal{x}^{- \delta} \big| \big|_{L^2 \rightarrow L^2} , \label{reductioncompacte}
\end{eqnarray}
since the norm in the right hand side goes to zero as the support of $ \Psi $ shrinks to $ \{ 0 \} $, for $ 0 $ is not an eigenvalue of $ \overline{P} $. The second order and first order term of (\ref{decompositioncompacte}), namely $ \overline{W}-V $, have a rather simple contribution. Indeed, we note that 
$$ || \scal{x}^M \big( \overline{W} - V \big) u ||_{L^2} \leq C || \overline{P}_0 u ||_{L^2}, \qquad u \in H^2, $$
using (\ref{1explicite}), (\ref{2explicite}) and (\ref{fort}) for $ P_0 $. By the Spectral Theorem $ \overline{P}_0 \big( \overline{P}_0 - i \epsilon \big)^{-1} $
is uniformly bounded on $ L^2 $ and thus
$$ || \scal{x}^M \big( \overline{W} - V \big)  \big( \overline{P}_0- i \epsilon \big)^{-1}  \Psi \big( \overline{P} \big) \scal{x}^{-\nu} ||_{L^2} \leq C || \Psi \big( \overline{P} \big) \scal{x}^{-\nu}  ||_{L^2 \rightarrow L^2} .  $$
We now consider $V$ alone. Since $ \scal{x}^M V $ has compact support, the Sobolev inequality and (\ref{faible}) for $ \overline{P}_0 $ yield
\begin{eqnarray}
  || \scal{x}^M V u ||_{L^2} \lesssim ||u||_{L^{2^*}} \lesssim || \nabla u ||_{L^2} \leq C \big| \big| \overline{P}_0^{1/2} u \big| \big|_{L^2 } . \label{star0}
\end{eqnarray}
We also observe that
\begin{eqnarray}
 || \overline{P}^{1/2} u ||_{L^2} \leq C || \overline{P}_0^{1/2} u ||_{L^2} , \label{star1}
\end{eqnarray}
since, by the compact support of $ \nabla \mbox{det} \ G (x) $  and the Sobolev inequality, 
$$ || \overline{P}^{1/2} u ||_{L^2}^2 \lesssim \big| \big| \nabla \big( \mbox{det} G(x)^{-1/4} u \big) \big| \big|_{L^2}^2
\lesssim  \left( || \nabla u ||_{L^2} + || u ||_{L^{2^*}} \right)^2 \lesssim || \nabla u ||^2_{L^2} . $$
Therefore, by (\ref{star0}) and (\ref{star1}), we have
\begin{eqnarray}
 || \scal{x}^M V  \big( \overline{P}_0 - i \epsilon \big)^{-1} \overline{P}^{1/2} u ||_{L^2} \leq C || u ||_{L^2}, \qquad \epsilon > 0 , \label{star2}
\end{eqnarray}
for all $u \in H^2 $. 
On the other hand, by approaching $ (\overline{P})^{-1/2} $ by
$$ S_n := \frac{1}{\sqrt{\pi}} \int_0^n e^{-tP} \frac{dt}{t^{1/2}}  , $$
in the sense that $ \overline{P}^{1/2} S_n \rightarrow I $ strongly on $ L^2 $, when applied to an $ H^2 $ function (see for instance \cite{Bouc1}), we deduce from (\ref{star2}) that
$$ || \scal{x}^M  V  \big( \overline{P}_0- i \epsilon \big)^{-1}  \Psi \big( \overline{P} \big) \scal{x}^{-\nu} ||_{L^2 \rightarrow L^2} \leq C \sup_{n} || S_n \Psi \big( \overline{P} \big)  \scal{x}^{-\nu}  ||_{L^2 \rightarrow L^2} . $$
Since $ S_n $ commutes with $ \Psi \big( \overline{P} \big) $, we shall obtain (\ref{reductioncompacte}) if we show that, for some $ \delta > 0 $,
$$ \sup_n || \scal{x}^{\delta} S_n \scal{x}^{-\nu} ||_{L^2 \rightarrow L^2} < \infty . $$
By the usual heat kernel bounds for $  \Delta_G $ ({\it e.g.} \cite{Davies} and references therein) and the fact that the Euclidean distance $ |x-y| $ is bounded from above and below by the geodesic distance $ d_G (x,y) $, we have
$$ [e^{-t\overline{P}}](x,y) \lesssim t^{-d/2} \exp (-|x-y|^2/Ct). $$ 
By integrating this estimate in $t$, we obtain that the kernel of $ \scal{x}^{\delta} S_n \scal{x}^{-\delta} $ satisfies
\begin{eqnarray*}
 0 \leq  \big[ \scal{x}^{\delta} S_n \scal{x}^{-\delta} \big](x,y) & \lesssim & \scal{x-y}^{\delta}|x-y|^{1-d} , \\
 & \lesssim &  |x-y|^{1-d + \delta} + f (x-y), 
\end{eqnarray*}
with $ f \in L^1 $. The convolution with $f$ is bounded on $ L^2 $ hence so is the operator with kernel $ f (x-y) \scal{y}^{\delta-\nu} $, if $ \delta \leq \nu $. We now consider the first term in the last line. By the Hardy-Littlewood-Sobolev inequality, the operator with kernel
$$ |x-y|^{1-d + \delta} \scal{y}^{\delta - \nu} , $$
is continuous on $ L^2 $ if $ \delta > 0 $ is small enough, since the convolution by $ |\cdot|^{1-d+\delta} $ maps $ L^{\frac{2d}{d+2+2\delta}} $ into $ L^2 $ and the multiplication by $ \scal{ \cdot}^{\delta-\nu} $ maps $ L^2 $ into $ L^{\frac{2d}{d+2+2\delta}} $. This shows that $ \scal{x}^{\delta}S_n \scal{x}^{-\delta + (\delta-\nu)} $ is uniformly bounded on $ L^2 $ and the result follows. \finpreuve
\bigskip

We consider now the second term of (\ref{Neumann2}).

\begin{prop} \label{absorbe2} Fix $ M > 0 $, $ \nu > 4 $ and $ \Psi \in C_0^{\infty} $. Then, if $ \epsilon_0 $ is small enough, we have
$$ \left| \left| \scal{x}^M \overline{W} \left( \big( \overline{P}_0-z \big)^{-1} - 
\big( \overline{P}_0- i \epsilon \big)^{-1} \right)  \Psi \big( \overline{P} \big) \scal{x}^{-\nu} \right| \right|_{L^2 \rightarrow L^2} \leq \frac{1}{4} , $$
for $ 0<|z|<\epsilon_0 $ (recall also the notation (\ref{notationz})).
\end{prop}

\noindent {\it Proof.} Recall first the standard fact that $ \Psi \big( \overline{P} \big) \scal{x}^{-\nu} $ preserves $ \scal{x}^{-\nu} $, ie $ \scal{x}^{\nu} \Psi \big( \overline{P} \big) \scal{x}^{-\nu} $ has a bounded closure on $ L^2$. Similarly $ \scal{x}^M \overline{W} (\overline{P_0}+1)^{-1} \scal{x}^{\nu} $ is bounded on $ L^2 $ since  $ \scal{x}^M \overline{W} $ has compact support and $ (\overline{P_0}+1)^{-1} $ preserves polynomial decay. It is thus sufficient to show that
$$ \left| \left| \scal{x}^{-\nu}  \big( \overline{P}_0 + 1 \big) \left( \big( \overline{P}_0-z \big)^{-1} - 
\big( \overline{P}_0- i \epsilon \big)^{-1} \right)   \scal{x}^{-\nu} \right| \right|_{L^2 \rightarrow L^2} \rightarrow 0, \qquad |z| \rightarrow 0 . $$
By writing
\begin{eqnarray*} 
 \big( \overline{P}_0 + 1 \big) \left( \big( \overline{P}_0-z \big)^{-1} - 
\big( \overline{P}_0- i \epsilon \big)^{-1} \right) & = &  \big( \overline{P}_0-z \big)^{-1} - 
\big( \overline{P}_0- i \epsilon \big)^{-1} + \frac{z}{\overline{P}_0 - z} - \frac{i \epsilon}{\overline{P_0}-i\epsilon} \\
& = & \int_0^{\lambda} (\overline{P}_0 - \mu - i \epsilon )^{-2} d \mu + \frac{z}{\overline{P}_0 - z} - \frac{i \epsilon}{\overline{P_0}-i\epsilon},
\end{eqnarray*}
the result follows from the bounds in Corollary \ref{corollairepratique1} for $ \overline{P}_0 $ with $ n = 1,2$, using in particular the integrability in $ \mu $ of 
$ ||\scal{x}^{-\nu}(\overline{P}_0 - \mu - i \epsilon )^{-2} \scal{x}^{-\nu}||_{L^2 \rightarrow L^2} $. \finpreuve

\bigskip
 
\noindent {\bf Proof of Theorem \ref{estimationprincipale}.} We first prove the theorem for some large enough $ \nu $ independent of $n$, namely $ \nu > 2N $. We shall see in the end of the proof how this implies the full result. So let us assume  that $ \nu > 2 N $. By Propositions \ref{absorbe1} and \ref{absorbe2}, by choosing $ \Psi $ with support close enough to $ 0 $ and by restricting $z$ to the region $ 0 < |z| \leq \epsilon_0 $ with $ \epsilon_0 $ small enough, we may assume that
 $$ \big| \big| \scal{x}^{\nu}  S_{\Psi} (z) \scal{x}^{-\nu} \big| \big|_{L^2 \rightarrow L^2} \leq 1/2 . $$
Therefore, by (\ref{pourNeumann}), we have
$$ || \scal{x}^{-\nu} R^{1}_{\psi} (z) \scal{x}^{-\nu} ||_{L^2 \rightarrow L^2} \leq \frac{4}{3} || \scal{x}^{-\nu} B^1 (z) \scal{x}^{-\nu} ||_{L^2 \rightarrow L^2} , $$
where we observe that the right hand side is bounded with respect to $z$: indeed, if we set more generally
$$ B^n (z) = \partial_z^{n-1} B^1 (z) , \qquad 1 \leq n \leq N , $$
 Corollary \ref{corollairepratique1} and Proposition \ref{estimeecotenegatif} for $ \overline{P}_0 $ show that
 we have
\begin{itemize}
\item{ if $ 1 \leq n \leq \bar{r}(d) $
\begin{eqnarray}
\big| \big| \scal{x}^{-\nu}  B^n (z) \scal{x}^{-\nu}   \big| \big|_{L^2 \rightarrow L^2} \leq C  \nonumber
\end{eqnarray}}
\item{if $ n = N = \bar{r}(d) + 1 $ and $ d $ is odd,
\begin{eqnarray}
\big| \big| \scal{x}^{-\nu}  B^{N} (z) \scal{x}^{- \nu}  \big| \big|_{L^2 \rightarrow L^2} \leq C |\mbox{Re}(z)|^{-1/2} , \nonumber  
\end{eqnarray}}
\item{if $ n= N $ and $d$ is  even, then for all $  2d < q < \infty $,  
\begin{eqnarray}
\big| \big| \scal{x}^{-\nu} B^{N}(z)  \scal{x}^{-\nu} \big| \big|_{L^{2} \rightarrow L^2} \leq C_{q} |\mbox{Re}(z)|^{- \frac{2d}{q}}  , \nonumber
\end{eqnarray}}
\end{itemize}
 for $ \mbox{Re}(z) \ne 0 $ and $ \mbox{Im}(z) \ne 0 $. We also have the same estimates for $ \scal{x}^{\nu} \partial_z^{n-1} S_{\Psi} (z) \scal{x}^{-\nu} $. In particular for $n=1$, this shows that
$$ || \scal{x}^{-\nu} R^{1}_{\psi} (z) \scal{x}^{-\nu} ||_{L^2 \rightarrow L^2} \lesssim 1 , \qquad 0 < |z| < \epsilon_0 . $$
For $n \geq 2$, we proceed as follows.  By applying $ \partial_z^{n-1} $ to (\ref{pourNeumann}), we obtain 
$$ \scal{x}^{-\nu} R_{\psi}^n(z) \scal{x}^{-\nu} = \scal{x}^{-\nu} \widetilde{B}^n (z) \scal{x}^{-\nu} +  \scal{x}^{-\nu}S_{\Psi}(\bar{z})^* R_{\psi}^n(z) S_{\Psi}(z) \scal{x}^{-\nu}  , $$
where, by an elementary induction, we see that $ \widetilde{B}^n (z) $ satisfy the same estimates  as $ B^n (z) $. Therefore
$$ || \scal{x}^{-\nu} R^n _{\psi} (z) \scal{x}^{-\nu} ||_{L^2 \rightarrow L^2} \leq \frac{4}{3} || \scal{x}^{-\nu} \widetilde{B}^n (z) \scal{x}^{-\nu} ||_{L^2 \rightarrow L^2} , $$
where the right hand side satisfies the expected estimates. We thus get the result with $ \nu > 2N $ for all $n = 1, \ldots, N$. To see that one can choose $ \nu > 2 n $, we proceed as follows. Fix $ M > 2 N $ and write, by  (\ref{pourNeumann}), 
$$ \scal{x}^{-\nu} R_{\psi}^1 (z) \scal{x}^{-\nu} = \scal{x}^{-\nu} B^1 (z) \scal{x}^{-\nu} + \scal{x}^{-\nu} S_{\Psi} (\bar{z})^* \scal{x}^M \left( \scal{x}^{-M} R_{\psi}^1 (z) \scal{x}^{-M} \right) \scal{x}^M S_{\Psi}(z) \scal{x}^{-\nu}. $$
We observe in this identity that $  \scal{x}^{-\nu} B^1 (z) \scal{x}^{-\nu} $ is bounded with respect to $z$, by the resolvent estimates for $ \overline{P}_0 $. The same holds for  $ \scal{x}^M S_{\Psi}(z) \scal{x}^{-\nu} $ since $ \scal{x}^M $ is harmless for $ \overline{W} $ has compactly supported coefficients. Therefore, the boundedness of $ \scal{x}^{-M} R_{\psi}^1 (z) \scal{x}^{-M} $ proved above gives the result for $n=1$. For $n \geq 2 $, we  differentiate  $ n-1 $ times with respect to $z$ and proceed as before. \finpreuve

\bigskip

\noindent {\bf Proof of Theorem \ref{estimationprincipaleLp}.} We may again replace $ (\overline{P}-z)^{-n} $ by its spectrally localized version (\ref{Resolvantepsi}) since $ (1 - \psi)(\overline{P}) (\overline{P}-z)^{-n} $ maps $ H^{-n} $ to $ H^n $, with bound independent of $z$ for small $z$,  and thus satisfies the expected $ L^{p} \rightarrow L^{p^{\prime}} $ boundedness. Let us consider first $n = 1$. Then $ p (1) = 2_* $ and $ q(1) = 2^* $. We start with (\ref{pourNeumann}) in which we observe that
\begin{eqnarray}
 || (\kappa A + i)^{-1} B^1 (z) (\kappa A - i )^{-1} ||_{L^{2_*} \rightarrow L^{2^*}} \lesssim 1 , \qquad 0 <|z| \leq 1 . \label{presqueconclusion1}
\end{eqnarray}
The estimate (\ref{presqueconclusion1}) follows from
\begin{eqnarray}
 || (\kappa A + i)^{-1} (\overline{P}_0 - z)^{-1} (\kappa A - i )^{-1} ||_{L^{2_*} \rightarrow L^{2^*}} \lesssim 1 , \qquad |\mbox{Re}(z)| \leq 1 , 
 \label{encorearappeler}
\end{eqnarray}
by Theorem \ref{theoremeprincipalpetit} and Proposition \ref{estimeecotenegatif} for $ P_0 $, and from 
\begin{eqnarray}
|| (\kappa A + i)^{-1} \Psi (\overline{P}) (\kappa A + i) ||_{L^{2^*} \rightarrow L^{2^*}} < \infty , \label{higherLebesgue1} \\
|| (\kappa A - i ) \overline{W} \psi (\overline{P}) (\kappa A + i) ||_{L^{2^*} \rightarrow L^{2_*}} < \infty , \label{higherLebesgue2}
\end{eqnarray}
as well as the adjoint estimates or similar ones with $ \psi $ instead of $ \Psi  $. The latter estimates follow easily from the fact that $ \Psi (\overline{P}) $ and $ \psi (\overline{P}) $ are pseudodifferential operators (see \cite{BTJFA}) with symbols in $ S^{-\infty}(\Ra^d \times \Ra^d) $ and the fact that $W $ has compactly supported coefficients.
We also have, for any $ \nu >2 $,
\begin{eqnarray}
 || \scal{x}^{\nu} \Psi (\overline{P}) S_{\Psi}(z) (\kappa A - i)^{-1} ||_{L^{2_*} \rightarrow L^{2}} \leq C , \qquad |\mbox{Re}(z)|<1 , \label{presqueconclusion2}
\end{eqnarray}
by (\ref{encorearappeler}) and the estimates
\begin{eqnarray}
 || \scal{x}^{\nu} \Psi (\overline{P}) \overline{W} (\kappa A + i) ||_{L^{2^*} \rightarrow L^2} < \infty , \label{higherLebesgue4} \\
|| (\kappa A - i) \Psi (\overline{P}) (\kappa A - i)^{-1} ||_{L^{2_*} \rightarrow L^{2_*}} < \infty , \label{higherLebesgue0}
\end{eqnarray}
which follow again from the fact that $ \Psi (\overline{P}) $ is a pseudodifferential operator of order $ - \infty $ and the compact support of the coefficients of $W$.
Therefore, by (\ref{pourNeumann}) where one can replace $ R_{\psi}^1 (z) $ by $ \Psi ( \overline{P} ) R^1_{\psi} (z) \Psi( \overline{P} ) $ in the right hand side,
and by (\ref{presqueconclusion1}) and (\ref{presqueconclusion2}), we obtain
$$ || (\kappa A + i)^{-1} R_{\psi}^1 (z) (\kappa A - i)^{-1} ||_{L^{2_*} \rightarrow L^{2^*}} \lesssim 1 + || \scal{x}^{-\nu} R_{\psi}^1 (z) \scal{x}^{-\nu} ||_{L^{2} \rightarrow L^{2}} , $$
so the result follows from Theorem \ref{estimationprincipale}. For $n \geq 2 $, we proceed by induction as in the proof of Theorem \ref{estimationprincipale} by applying $ \partial_z^{n-1} $ to (\ref{pourNeumann}). We omit the details but rather point out that the analogues of the estimates (\ref{higherLebesgue1}), (\ref{higherLebesgue2}), (\ref{higherLebesgue4}) and (\ref{higherLebesgue0}) associated to $ q (n) $ don't cause any trouble when $ q (n) = q (N) = \infty $ since they involve pseudodifferential operators of order $ - \infty $ (but no zero order pseudodifferential operator) which are bounded on all $ L^p $ spaces for $ p \in [1,\infty] $. \finpreuve

\section{Local energy decay} \label{referenceStone}
\setcounter{equation}{0}
The purpose of this section is to prove Theorems \ref{theoremeexponentiel} and \ref{theoremenontrapping}. For convenience, we work with the self-adjoint realization $ \overline{P} $ on $L^2 (\Ra^d,dx) $ of
$$  P = - \mbox{det}\ G(x)^{-1/4} \frac{\partial}{\partial x_j} \left( \mbox{det}\ G(x)^{1/2}  G_{jk}(x)  \frac{\partial}{\partial x_k} \right) \mbox{det} \ G (x)^{-1/4} ,  $$
which is unitarily equivalent to $ - \Delta_G $ on $ L^2 (\Ra^d,d_G x) $.

\subsection{Spectral localization}
Let $ m \geq 0 $ be a real number and $ \alpha = 1 $ or $ 1/2 $. In this paragraph, we define  $ U (t) $ by 
$$ U (t) = e^{it (\overline{P}+m^2)^{\alpha}}  , $$
which will allow to cover simultaneously  the Schr\"odinger ($m = 0$, $ \alpha = 1 $), Klein-Gordon ($ m > 0 $, $ \alpha = 1/2 $) and wave equations ($ m =0 $, $ \alpha = 1/2 $). Our purpose here is to reduce estimates on such  flows to spectrally localized estimates. Actually, the result of this subsection only uses that $ U (t) $ is some bounded function of $ \overline{P} $ and nothing else. 

Consider a dyadic partition of unit
\begin{eqnarray}
 1  &= &  \Phi_0 (\lambda) + \sum_{k \geq 0} \varphi (2^{-k} \lambda) \\
  & = & \Phi_0 (\lambda) + \Phi (\lambda), \label{hautefrequences}
\end{eqnarray}
defined for $ \lambda $ near $ [0 , \infty ) $, with
$$ \Phi_0 \in C_0^{\infty} (\Ra), \qquad \varphi \in C_0^{\infty}(0,+\infty) . $$
We also select $ \psi $ such that
\begin{eqnarray}
 \psi \in C_0^{\infty} (0,+ \infty), \qquad \psi \equiv 1 \ \ \mbox{near} \ \mbox{supp}(\varphi) . \label{definitionpsi}
\end{eqnarray}
It will be convenient to denote
$$ E_{\nu} (h,t) =  \langle x \rangle^{-\nu} U(t) \varphi (h^2 \overline{P}) \langle x \rangle^{-\nu} , $$
and
\begin{eqnarray}
 e_{\nu} (h,t) = \big| \big| E_{\nu}(h,t)  \big|  \big|_{L^2 \rightarrow L^2} . \label{notationnorme}
\end{eqnarray}
Our main purpose here is to show the following proposition.
\begin{prop} \label{localisationspectrale} For all $ \nu \geq 0 $ and $ M > 0  $ there exists $ C > 0 $ such that
\begin{eqnarray*}
 \big| \big|  \langle x \rangle^{-\nu} \Phi(\overline{P}) U (t) \langle x \rangle^{-\nu}  u \big| \big|_{ L^2}^2  \leq C  \sum_{h^2 = 2^{-k}} 
  e_{\nu}(h,t)^2 \left( \big| \big| \psi (h^2 \overline{P}) u \big| \big|^2_{L^2} + h^{M} \big| \big| (1-\overline{P})^{-M/2} u \big| \big|_{L^2}^2 \right)  ,
\end{eqnarray*}
for all $ t \in \Ra $ and  $ u \in {\mathcal S}(\Ra^d) $. Here $ \Phi $ is defined in (\ref{hautefrequences}).
\end{prop}
As a corollary, we obtain the following estimate which we shall use in Subsection \ref{decroissancetemporelle}.
\begin{coro} \label{corollairefinal} For all $ \nu \geq 0 $ and $ s \in \Ra $, one has
$$ \big| \big| \langle x \rangle^{-\nu} \Phi(\overline{P}) U(t) \langle x \rangle^{-\nu} u \big| \big|_{ L^2} \leq C_{\nu,s} \left( \sup_{h \in (0,1]}  h^s e_{\nu}(h,t)  \right) || u ||_{H^s} , $$
for all $ t \geq 0 $ and $ u \in {\mathcal S} (\Ra^d) $.
\end{coro}

\noindent {\it Proof of Corollary \ref{corollairefinal}.} By the Spectral Theorem, we have 
$$  h^{-s} || \psi (h^2 \overline{P}) u  ||_{L^2} \leq C || \psi (h^2 \overline{P}) (1-\overline{P})^{s/2} u ||_{L^2} , $$
for all $ u \in L^2 $. Since $ || (1-\overline{P})^{s/2} u ||_{L^2} \leq C || u ||_{H^s} $ by classical elliptic estimates, we obtain, by almost orthogonality,
\begin{eqnarray}
 \sum_{h^2 = 2^{-k}} h^{-2s} || \psi (h^2 \overline{P}) u ||_{L^2}^2 \leq C || u ||_{H^s}^2 . \label{sommequasiorthogonale}
\end{eqnarray}
On the other hand, by Proposition \ref{localisationspectrale}, we have
$$ \big| \big| \langle x \rangle^{-\nu} \Phi(\overline{P}) U(t) \langle x \rangle^{-\nu} u \big| \big|_{ L^2}^2 \lesssim  \sum_{h^2 = 2^{-k}} 
  e_{\nu}(h,t)^2 h^{2s} \left( h^{-2s} \big| \big| \psi (h^2 \overline{P})  u \big| \big|^2_{L^2} + h^{M-2s} \big| \big|  u \big| \big|_{H^{-M}}^2 \right) .  $$
Choosing $ M > 2|s| $, we have $ \sum_h h^{M-2s} < \infty $, $ || u ||_{H^{-M}} \leq || u ||_{H^s} $ and we conclude using (\ref{sommequasiorthogonale}). 
 \finpreuve

\bigskip

We now consider the proof of Proposition \ref{localisationspectrale}. Write first
$$ \langle x \rangle^{-\nu} \Phi(\overline{P}) U(t) \langle x \rangle^{-\nu} u = \sum_{h^2 = 2^{-k}} E_{\nu} (h,t) u , $$
where the sum converges weakly (and actually in $ L^2 $ by the analysis below). We will need the following  result.
\begin{lemm} \label{lemmepratique} For all $ M \geq 0 $, one has
$$ \varphi (h^2 \overline{P}) \langle x \rangle^{-\nu} (1-\psi(h^2\overline{P})) = h^M \varphi (h^2 \overline{P}) \scal{x}^{-\nu} R_{M,\nu} (h) , $$
with
$$ \big| \big| R_{M,\nu}(h) (1-\overline{P})^{M/2} \big| \big|_{L^2 \rightarrow L^2} \leq C, \qquad h \in (0,1 ] . $$
\end{lemm}

\noindent {\it Proof.} By (\ref{definitionpsi}), we can select $ \widetilde{\varphi} \in C_0^{\infty} (0,+\infty) $ such that
$$ \varphi \widetilde{\varphi} = \varphi  \qquad \mbox{and} \qquad \psi \equiv 1 \ \ \mbox{near} \ \mbox{supp}(\widetilde{\varphi}) , $$
 and thus write
$$ \varphi (h^2 \overline{P}) \langle x \rangle^{-\nu} = \varphi (h^2 \overline{P}) \langle x \rangle^{-\nu} \big(\scal{x}^{\nu} \widetilde{\varphi}(h^2\overline{P}) \scal{x}^{-\nu} \big) . $$
The result follows then from the fact that, for all $ M $,
$$ \big| \big| \big(\scal{x}^{\nu} \widetilde{\varphi}(h^2\overline{P}) \scal{x}^{-\nu} \big) (1 -  \psi(h^2 \overline{P})) (1-\overline{P})^{M/2} \big| \big|_{L^2 \rightarrow L^2} \leq C_{M,\nu} h^M , \qquad h \in (0,1 ] , $$
by pseudodifferential functional calculus ({\it e.g.} \cite{BTJFA}), since all terms of the pseudo-differential expansion cancel because $ \widetilde{\varphi} $ and $ 1 - \psi $ have disjoint supports. \finpreuve

\bigskip

\noindent {\bf Proof of Proposition \ref{localisationspectrale}.}  By Lemma \ref{lemmepratique}, we have
$$ E_{\nu}(h,t) = \psi (h^2 \overline{P}) E_{\nu}(h,t) \psi (h^2 \overline{P})  + h^M R_{M,\nu}(h)^* E_{\nu}(h,t) \psi (h^2 \overline{P})  + h^M E_{\nu}(h,t) R_{M,\nu}(h) , $$
and the result will follow from the estimates on each term given below.
\begin{enumerate} 
\item{{\it 1st term.} By almost orthogonality, we have
\begin{eqnarray*}
 \left| \left| \sum_{h^2 = 2^{-k}} \psi (h^2 \overline{P}) E_{\nu}(h,t) \psi (h^2 \overline{P}) u \right| \right|_{L^2}^2 & \lesssim & \sum_{h^2 = 2^{-k}} || E_{\nu}(h,t) \psi (h^2 \overline{P}) u ||_{L^2}^2 , \\
 & \lesssim & \sum_{h^2 = 2^{-k}}  e_{\nu}(h,t)^2 ||\psi (h^2 \overline{P}) u ||_{L^2}^2 .
\end{eqnarray*}}
\item{{\it 2nd term.} Since $ || R_{M,\nu} (h) ||^*_{L^2 \rightarrow L^2} \leq C $ by Lemma \ref{lemmepratique}, we also have
\begin{eqnarray*}
\big| \big| \sum_{h^2 = 2^{-k}} h^M R_{M,\nu}(h)^* E_{\nu}(h,t) \psi (h^2 \overline{P}) u \big| \big|_{L^2} & \lesssim & \sum_{h^2 = 2^{-k}} h^M \big| \big| E_{\nu}(h,t) \psi (h^2 \overline{P}) u \big| \big|_{L^2} , \\ & \lesssim & \left(
 \sum_{h^2 = 2^{-k}} e_{\nu}(h,t)^2  \big| \big| \psi (h^2 \overline{P}) u \big| \big|_{L^2}^2 \right)^{1/2} ,
\end{eqnarray*}
by the Cauchy-Schwarz inequality since $ \sum_{h^{2} = 2^{-k}} h^{2M} < \infty $.}
\item{{\it 3rd term}. By Lemma \ref{lemmepratique}, 
\begin{eqnarray*}
\big| \big| \sum_{h^2 = 2^{-k}} h^M E_{\nu}(h,t) R_{M,\nu}(h) u \big| \big| & \lesssim & \sum_{h^2 = 2^{-k}} h^M e_{\nu}(h,t) || (1-\overline{P})^{-M/2} u ||_{L^2} , \\ & \lesssim & \left( \sum_{h^2 = 2^{-k}} h^M e_{\nu}(h,t)^2 || (1-\overline{P})^{-M/2} u ||_{L^2}^2 \right)^{1/2} ,
\end{eqnarray*}
again by the Cauchy-Schwartz inequality since $ \sum_{h^2 = 2^{-k}} h^M < \infty $.}
\end{enumerate}
The proof is complete.  \finpreuve

\subsection{Semiclassical estimates} \label{resolvanteverspuissances}
To prove quantitative decay rates for the Schr\"odinger group, we shall use integration by parts in the Stone formula. For this purpose, we need to estimate powers of the resolvent. In this subsection, we show that, if one has semiclassical estimates for the resolvent, then one has estimates for its powers. For simplicity, we will only consider the square of the resolvent, but higher powers  can be treated similarly.

We introduce the usual notation
$$ R (z,h) = (h^2 \overline{P}-z)^{-1} . $$
Throughout this subsection, $ J_0 \Subset (0,\infty) $ will be a relatively compact interval satisfying the following condition.

\noindent {\bf Assumption A}. {\it There exist a real number $ \nu_0 \geq 0 $ and a function $ F : (0,1] \rightarrow (0,+\infty) $ satisfying
\begin{eqnarray}
 F (h) \gtrsim h^{-1} , \label{borneinfoptimale}
\end{eqnarray}
such that, for all $ \nu > \nu_0 $ and all open interval $ J \Subset J_0 $,}
 \begin{eqnarray}
 || \langle x \rangle^{- \nu} R(z,h) \langle x \rangle^{-\nu} ||_{L^2 \rightarrow L^2} \leq C_{\nu,J} F (h), \qquad h \in (0,1] ,  \ \mbox{Re}(z) \in J . \label{exponentielle2}
\end{eqnarray}

\bigskip

Without any condition on $ G $, such estimates holds with $ F (h) = C e^{C/h} $ (\cite{Burq1,Burq2} and \cite{stCardosoVodev}). When the geodesic flow is non trapping, one can choose $ F (h) = C / h $ \cite{RoTa,Wang,GeMa,RoENS,VaZw}. In some cases where one has weak trapping one may take $ F (h) = C |\log h|/h $ or polynomial powers of $ h^{-1} $ \cite{Naka1,NoZw}.

Our purpose here is to prove the following.
\begin{prop} \label{deriveeexponentielle} If Assumption A holds then, for  all $ \nu > \nu_0 $ and all interval $ J \Subset J_0 $, there exists $ C > 0 $ such that
\begin{eqnarray}
 || \langle x \rangle^{- \nu-1} R (z,h)^2 \langle x \rangle^{-\nu-1} ||_{L^2 \rightarrow L^2} \leq C F (h)^2 , \label{exponentielle2bis}
\end{eqnarray}
for all $h \in (0, 1]$ and all $ z $ such that $ \emph{Re}(z) \in J $. 
\end{prop} 


\bigskip


The principle of the proof below is well known (see \cite{IKRem} and \cite{AJensen}) but we recall  the main steps to emphasize the behaviour with respect to $ h $ (the previous works addressed either the case $h=1$ or the high energy limit for potentials, which is a non trapping case).
The approach is based on microlocal parametrices of the semiclassical Schr\"odinger group $ e^{-ith\overline{P}} $, from which we recover the resolvent by
\begin{eqnarray}
R(z,h) = \frac{i}{h} \int_0^{\pm \infty} e^{it z/h} e^{-ith\overline{P}} dt, \qquad \pm \mbox{Im}(z) > 0 . \label{resolventepropagateur}
\end{eqnarray}
It is  convenient to record the following elementary lemma.
\begin{lemm} \label{Fubinistrong} Let $ A (t) , B (t) $ be bounded operators on $ L^2 (\Ra^d) $, strongly continuous with respect to  $t$  and such that, for some $ N \geq 0$,
$$ || A (t) ||_{L^2 \rightarrow L^2} + || B (t) ||_{L^2 \rightarrow L^2} \leq C \langle t \rangle^N, \qquad t \in \Ra . $$
Then
\begin{eqnarray}
\int_0^{\pm \infty} e^{it\zeta} \left( \int_0^t A (t-s)B(s) ds \right) dt = \left( \int_0^{\pm \infty} e^{it\zeta} A (t)dt \right) \left( \int_0^{\pm \infty} e^{it\zeta} B (t)dt \right) , \nonumber
\end{eqnarray}
provided that
$$ \pm \emph{Im}(\zeta) > 0 . $$
\end{lemm}

We will use the well known Isozaki-Kitada parametrix, introduced first for potential scattering (see \cite{IsKi}). Here we need it in the metric case with a semiclassical parameter. In this context, we refer for instance to \cite{BTJFA} for the details or proofs of the statements quoted below, in particular Lemma \ref{IsozakiKitada}. We recall only what is necessary for the proof of Proposition \ref{deriveeexponentielle}.

Denote by $ S_{\rm scat} (\mu,-\infty) $  the set of smooth functions $ a $ on $ \Ra^{2d} $ such that, for all $ M > 0 $,
$$ \big| \partial_x^{\alpha} \partial_{\xi}^{\beta} a (x,\xi) \big| \leq C_{\alpha \beta M} \scal{x}^{\mu-|\alpha|} \scal{\xi}^{-M} ,  $$
where the best constants $ C_{\alpha \beta M} $ are seminorms for which it is a Fr\'echet space. 

Given real numbers $ R > 0 $, $ \sigma \in ( - 1 , 1 ) $ and any interval $ I \Subset (0,+\infty) $, one defines the outgoing ($+$) and incoming ($-$) areas by
$$ \Gamma^{\pm} (R,I,\sigma) := \big\{ (x,\xi) \in \Ra^{2d} \ | \ |x|>R, \ |\xi|^2 \in I, \  \pm x \cdot \xi > \sigma |x| |\xi| \big\} . $$
It turns out that, for any $ I $ and $ \sigma $ as above, one can choose $ R $ large enough so that one can solve the following eikonal equations
$$ \nabla_x \varphi^{\pm}(x,\xi) \cdot G (x)^{-1} \nabla_x \varphi^{\pm} (x,\xi) = |\xi|^2 , $$
for $ (x,\xi) \in \Gamma^{\pm}(R,I,\sigma) $ with solutions which are close to the free phase $ x \cdot \xi $ ({\it i.e.} the solution if $ G \equiv I $) in the sense that
$$ \big| \partial_x^{\alpha} \partial_{\xi}^{\beta} ( \varphi^{\pm}  (x,\xi) - x \cdot \xi ) \big| \leq C_{\alpha \beta} \scal{x}^{1-\rho-|\alpha|}, \qquad (x,\xi) \in \Gamma^{\pm}(R,I,\sigma) , $$
where $ \rho > 0 $ is the same as in (\ref{cometriquelongueportee}). One can then define the following Fourier integral operators
$$ J_{\pm}(a^{\pm})u(x,h) = (2 \pi h)^{-d} \int \! \! \int e^{\frac{i}{h}( \varphi^{\pm}(x,\xi)-y \cdot \xi )} a^{\pm}(x,\xi) u (y) dy d \xi , $$
for symbols such that
$$ a^{\pm} \in S_{\rm scat} (0,-\infty), \qquad \mbox{supp}(a^{\pm}) \in \Gamma^{\pm}(R,I,\sigma) . $$
We can now give a form of the Isozaki-Kitada parametrix.
\begin{lemm}[Isozaki-Kitada parametrix] \label{IsozakiKitada} Fix two intervals $ I \Subset I^{\prime} \Subset (0,+\infty) $. Then, for all $ R $ large enough and all
$$ \chi^{\pm} \in S_{\rm scat}(0,-\infty), \qquad \emph{supp}(\chi^{\pm}) \subset \Gamma^{\pm}(R,I,-1/2), $$
we can find, for all $ M \geq 0 $, symbols
\begin{eqnarray*}
 a^{\pm}_M (h) & \in & S_{\rm scat}(0,-\infty),  \qquad  \emph{supp}(a^{\pm}_M(h)) \subset \Gamma^{\pm}\big( R^{1/4},I^{\prime},-9/10 \big), \\
 b^{\pm}_M (h) & \in & S_{\rm scat}(0,-\infty),  \qquad  \emph{supp}(b^{\pm}_M(h)) \subset \Gamma^{\pm}\big( R^{1/2},I^{\prime},-3/4 \big), \\
 r^{\pm}_M (h) & \in & S_{\rm scat}(-2M,-\infty) , 
\end{eqnarray*} 
bounded with respect to $h$ in their classes, such that
\begin{enumerate}
\item{\begin{eqnarray*}
 e^{-ith\overline{P}} \chi^{+}(x,hD) & = & J_{+}\big( a^{+}_M(h) \big) e^{ith\Delta} J_{+} \big( b^{+}_M(h) \big)^* + h^M e^{-ith\overline{P}} r_{M}^{+}(x,hD,h) \\
 & & + h^{M} \int_0^t e^{-i(t-s)h\overline{P}} B_{M}^{+}(s,h) ds ,
\end{eqnarray*}
with $ B_M^{+} (s,h) $ strongly continuous with respect to $s$ and such that
$$ \big| \big| \langle x \rangle^{M} B_M^{+} (s,h) \langle x \rangle^{M}  \big| \big|_{L^2 \rightarrow L^2} \leq C \langle s \rangle^{-M}, \qquad \ s \geq 0 , \ h \in (0,1] . $$}
\item{(Adjoint case)
\begin{eqnarray*}
 \chi^{-} (x,hD) e^{-ith\overline{P}} & = & J_{-}\big( b^{-}_M(h) \big) e^{ith\Delta} J_{-} \big( a^{-}_M(h) \big)^* + h^M  r_{M}^{-}(x,hD,h) e^{-ith\overline{P}} \\
 & & + h^{M} \int_0^t B_{M}^{-}(-s,h) e^{-i(t-s)h\overline{P}} ds ,
\end{eqnarray*}
with $ B_M^{-} (-s,h) $ strongly continuous with respect to $s$ and such that
$$ \big| \big| \langle x \rangle^{M} B_M^{-} (-s,h) \langle x \rangle^{M}  \big| \big|_{L^2 \rightarrow L^2} \leq C \langle s \rangle^{-M}, \qquad  s \geq 0 , \ h \in (0,1] . $$}
\end{enumerate}
\end{lemm}
We simply point out that this lemma gives  good approximations for $ t \geq 0 $ only, which will be sufficient for us. There is of course a similar statement for negative times by exchanging $ + $ and $ - $ everywhere. 

We also mention that the symbols $ a^{\pm}_M(h) $ and $ b^{\pm}_M (h) $ are finite sums of the form $ \sum h^j c^{\pm}_j $ with $ c_j \in S_{\rm scat} (-j,-\infty) $ independent of $h$. The following lemma will thus be useful to estimate the leading terms of the parametrix. Again, we consider only positive times.

\begin{lemm}[Free propagation estimates] \label{freeprop} Let $ \mu_1 \geq \mu_2 \geq 0 $ be real numbers, $ I \Subset (0,+\infty) $ an interval and $ \sigma \in (-1,1) $. Then, for all $ R  $ large enough and all symbol $ c^{\pm} $ satisfying,
$$ c^{\pm} \in S_{\rm scat}(0,-\infty), \qquad \emph{supp}(c^{\pm}) \subset \Gamma^{\pm}(R,I,\sigma) , $$
 we have
\begin{enumerate}
\item{
$$
 \big| \big| \langle x \rangle^{-\mu_1} e^{ith\Delta}  J_+ (c^+)^* \langle x \rangle^{\mu_2} \big| \big| \leq C \langle t \rangle^{\mu_2 - \mu_1}, \qquad  t \geq 0 , \ \ h \in (0,1] , $$}
\item{(Adjoint case)
$$ \big| \big| \langle x \rangle^{\mu_2} J_- (c^-) e^{ith\Delta}  \langle x \rangle^{-\mu_1} \big| \big| \leq C \langle t \rangle^{\mu_2 - \mu_1}, \qquad  t \geq 0 , \ \ h \in (0,1] . $$}
\end{enumerate}
\end{lemm}

We refer for instance to \cite{IKRem} or \cite{BTJFA} for a proof of this lemma, which is fairly elementary and follows from integrations by parts in the (explicit) kernel of the operators for integers $ \mu_1,\mu_2 $ and then by an interpolation argument for real ones.

\bigskip

\noindent {\bf Proof of Proposition \ref{deriveeexponentielle}.} We may assume that $ \mbox{Im}(z) > 0 $, otherwise one takes the adjoint. By the Spectral Theorem, it is  sufficient to prove a $ {\mathcal O}(F(h)^2) $ upper bound for
$$ \phi (h^2 \overline{P}) R (z,h)^2 = R (z,h) \phi (h^2\overline{P}) R (z,h) , $$
with $ \phi \in C_0^{\infty}(0,+\infty) $ which is equal to $1$ near the interval $ J $ where $ \mbox{Re}(z) $ lives.
  Let $ \chi \in C_0^{\infty} (\Ra^d) $ such that $ \chi (x) = 1 $ for $ |x| \leq R $, with $R$ to be chosen below according to Lemmas \ref{IsozakiKitada} and \ref{freeprop}. Then
\begin{eqnarray}
 \phi (h^2 \overline{P}) R (z,h)^2 =  R (z,h)  \phi (h^2 \overline{P})\chi R (z,h) +  R (z,h) \phi (h^2 \overline{P}) (1-\chi) R (z,h) . \label{decoupageresolvente1}
\end{eqnarray}
Since, for any $ M > 0 $,
$$ \big| \big| \scal{x}^M \phi (h^2 \overline{P})\chi \scal{x}^M \big| \big|_{L^2 \rightarrow L^2} \lesssim 1, \qquad h \in (0,1] , $$
the $ F (h)^2 $ upper bound for the first term in the right hand side of (\ref{decoupageresolvente1}), weighted on both sides by $ \scal{x}^{-\nu-1} $, follows easily from (\ref{exponentielle2}). Note that the extra power $ \scal{x}^{-1} $ is useless for this term. In the second term, we use the following pseudodifferential expansion (see \cite{BTJFA}): for all $ M \geq 1 $,
$$ \phi (h^2 \overline{P}) (1-\chi) = \sum_{j < M} h^j \chi_{j}^{+}(x,hD) + \sum_{j< M} \chi_{j}^{-}(x,hD) + h^M R_M (h) , $$
where, if $ I \Subset (0,+\infty) $ is a neighborhood of $ \mbox{supp}(\phi) $ and $R$ is large enough,
$$ \chi_{j}^{\pm} \in S_{\rm scat}(-j,-\infty) , \qquad \mbox{supp}\big(\chi_{j}^{\pm} \big) \in \Gamma^{\pm} (R,I,-1/2) , $$
and
$$ \big| \big| \scal{x}^{M/2} R_M (h) \scal{x}^{M/2} \big| \big|_{L^2 \rightarrow L^2} \lesssim 1, \qquad h \in (0,1] . $$
By choosing $ M $ large enough, the contribution of $ R_M (h) $ is treated similarly to the one of $ \phi (h^2 \overline{P}) \chi $ above, so we are left with the study of terms of the form
$$ \scal{x}^{-\nu-1} R (z,h) \chi^{\pm} (x,hD) R (z,h) \scal{x}^{-\nu-1} . $$
 The idea is to use Lemma \ref{IsozakiKitada} for
$$ R (z,h) \chi^{+}(x,hD) \qquad \mbox{and} \qquad \chi^{-}(x,hD) R (z,h) , $$
by expanding $ R (z,h) $ via (\ref{resolventepropagateur}), with  $ t \geq 0$ since $ \mbox{Im}(z)>0 $. We consider $ \chi^+ $.
By Lemma \ref{Fubinistrong} and item 1 of Lemma \ref{IsozakiKitada}, we have
$$ R (z,h) \chi^{+}(x,hD) = J_{+}\big( a^{+}_M(h) \big) (-h^2\Delta-z)^{-1} J_{+} \big( b^{+}_M(h) \big)^* + h^{M-1} R (z,h) R^+_M (h), $$
where
$$ R^+_M (h) = h r_M^+ (x,hD,h) + \int_0^{+ \infty} e^{-itz/h} B_{M}^{+}(t,h) dt , $$
satisfies
$$ \big| \big| \scal{x}^{M} R^+_M (h) \scal{x}^{M/8} \big| \big|_{L^2 \rightarrow L^2} \leq C, \qquad h \in (0,1] . $$
The contribution of $ R^{+}_M (h) $ is thus similar to the one of $ R_M (h) $ and $ \phi (h^2 \overline{P}) \chi $ above. We then consider
$$ \scal{x}^{-\nu-1} \left( J_{+}\big( a^{+}_M(h) \big) (-h^2\Delta-z)^{-1} J_{+} \big( b^{+}_M(h) \big)^* \right) R (h,z) \scal{x}^{-\nu - 1} . $$
Choose $ \nu^{\prime} $ such that
$$ \nu_0 < \nu^{\prime} < \nu . $$
Then, by Assumption A,
$$ \big| \big| \scal{x}^{-\nu^{\prime}} R (h,z) \scal{x}^{-\nu - 1} \big| \big| \leq C F (h) , \qquad h \in (0,1] . $$ On the other hand, using item 1 of Lemma \ref{freeprop}, (\ref{resolventepropagateur}) for $ - \Delta $ and (\ref{borneinfoptimale}), we have
\begin{eqnarray*}
 \left| \left| \scal{x}^{-\nu-1} \left( J_{+}\big( a^{+}_M(h) \big) (-h^2\Delta-z)^{-1} J_{+} \big( b^{+}_M(h) \big)^* \right) \scal{x}^{\nu^{\prime}} \right| \right|_{L^2 \rightarrow L^2} & \lesssim &  h^{-1} , \\
 & \lesssim & F (h)  .
\end{eqnarray*}
Here we use the additional fact that $ \scal{x}^{-\nu-1} J_+ (a^+_M (h)) \scal{x}^{\nu +1} $ is bounded on $ L^2 $, uniformly in $h$.
All this shows that $ || \scal{x}^{-\nu-1} R (z,h) \chi^+ (x,hD) R (z,h) \scal{x}^{-\nu-1} ||_{L^2 \rightarrow L^2} $ is bounded by $ C F (h)^2 $. The same analysis holds for $ \chi^- $ using the Adjoint Cases in Lemma \ref{IsozakiKitada} and  Lemma \ref{freeprop} and this completes the proof. \finpreuve

\subsection{Time decay} \label{decroissancetemporelle}

In this paragraph, we prove Theorem \ref{theoremeexponentiel} and Theorem \ref{theoremenontrapping}.




The following proposition will give the contribution of the low frequencies.

\begin{prop} \label{etape1prop} Let $ m > 0 $ and $ \chi \in C_0^{\infty} (\Ra, \Ra) $. For each $ t \in \Ra $, let $ \varphi_t (\lambda) $ denote any of the following functions
$$  \chi (\lambda) e^{-it \lambda},  \qquad \chi (\lambda) e^{-it (|\lambda|+m^2)^{1/2}}, \qquad \chi (\lambda ) \cos \left( t |\lambda|^{1/2} \right) , \qquad \chi (\lambda) \frac{\sin \left( t |\lambda|^{1/2} \right) }{|\lambda|^{1/2}} .  $$
Then, for all $ \nu > 2 (\bar{r}(d)+1) $, there exists $ C $ such that
\begin{eqnarray}
 \big| \big| \scal{x}^{-\nu} \varphi_t (\overline{P}) \scal{x}^{-\nu} \big| \big|_{L^2 \rightarrow L^2} \leq C \scal{t}^{-\bar{r}(d)}, \qquad t \in \Ra . \label{borneainterpoler} 
\end{eqnarray}
\end{prop}

\bigskip

\noindent {\bf Remark.} We note that the $ L^2 \rightarrow L^2 $ estimate of this proposition can be turned into a $  H^{-s} \rightarrow H^{s} $ estimate for all $ s \geq 0 $. Indeed, one can write
$$ \varphi_t (\overline{P}) = \widetilde{\chi} (\overline{P})\varphi_t (\overline{P}) \widetilde{\chi} (\overline{P}) , $$
with $ \widetilde{\chi} \in C_0^{\infty} $ such that $ \widetilde{\chi} \chi = \chi $, and use the fact that for any $ \nu \geq 0 $
$$ \scal{x}^{-\nu} \widetilde{\chi} (\overline{P}) \scal{x}^{\nu} : L^2 \rightarrow H^s , $$
is bounded (which follows for instance from the form of $ \widetilde{\chi}(\overline{P}) $ given in \cite{BTJFA}).

\bigskip

We quote the following result whose proof can be found in \cite{ReedSimon}.
\begin{lemm}[Stone's formula] \label{Stones} For all compactly supported continuous function $ \varphi \in C_0^0 (\Ra) $, one has
$$ \varphi \big( \overline{P} \big) = \lim_{\delta \downarrow 0 } \frac{1}{\pi} \int_{\Ra} \varphi (\lambda) \ \emph{Im} (\overline{P}-\lambda - i \delta)^{-1} d\lambda, $$
the limit being taken in the strong sense. Here $ \emph{Im}B = \frac{B-B^*}{2i} $.
\end{lemm}

\bigskip

\noindent {\it Proof of Proposition \ref{etape1prop}.} Denote for simplicity $ r = \bar{r}(d) $. In the first case,  $r$ integrations by part in the integral  yield
\begin{eqnarray}
 t^{r} \int_{\Ra} \varphi_t (\lambda) \mbox{Im} (\overline{P}- \lambda - i \delta)^{-1} d \lambda  = (-i)^r \int_{\Ra} e^{-it \lambda} \partial^r_{\lambda} \left( \chi (\lambda) \mbox{Im} (\overline{P}- \lambda - i \delta)^{-1} \right) d \lambda . \label{explicationIPP}
\end{eqnarray}
By Theorem \ref{estimationprincipale}, the  first $r$ derivatives with respect to $ \lambda $ of $(\overline{P}- \lambda \pm i \delta)^{-1} $ are integrable near $ 0 $, in the suitable weighted spaces, with uniform bounds in $ \delta $. Thus the right hand side of (\ref{explicationIPP}) is bounded uniformly with respect to $ \delta $ and $t$ and the result follows by using Lemma \ref{Stones}. The second case is similar, once we have noticed the following points. Since $ \overline{P} $ is non negative, we may modify $ \chi $ as we wish on $ (- \infty,0) $ without changing the operator $ \varphi_t (\overline{P}) $. In particular, we may assume that $ \chi $ is supported in $ \{ \lambda  > - m^2/2 \} $ and then
$$ \varphi_t (\overline{P}) = \chi (\overline{P}) e^{it (\overline{P}+m^2)^{1/2}} = \widetilde{\varphi}_t (\overline{P}) , $$
with $ \widetilde{\varphi}_t (\lambda) =  \chi (\lambda) e^{-it (\lambda+m^2)^{1/2}} $. The result follows again by integrating by part, using
$$ t  e^{-it (\lambda+m^2)^{1/2}} =2 i(\lambda+m^2)^{1/2} \partial_{\lambda} e^{-it (\lambda + m^2)^{1/2}} $$
on the support of $ \chi $ where $ \lambda + m^2 > m^2/2 $. In the last two cases, we need to work a little bit more since we shall have boundary terms in the integrations by part. We treat the last case, the third one being similar. By setting 
$$ B_{\delta} (\lambda) = \chi (\lambda) (\overline{P}-\lambda-i\delta)^{-1}, $$
and by  the change of variables $ \lambda = \pm \mu^2 $ on $ \Ra^{\pm} $, we have
\begin{eqnarray}
 \int_{\Ra} \frac{\sin (t|\lambda|^{1/2})}{|\lambda|^{1/2}} \chi (\lambda) \mbox{Im} (\overline{P}- \lambda - i \delta)^{-1} d \lambda = 2 \int_{0}^{\infty}  \sin (t \mu)  \mbox{Im} \left( B_{ \delta} (\mu^2) + B_{\delta} (-\mu^2) \right)  d \mu . \label{righthandsideIPP}
 \end{eqnarray}
By $r$ integrations by part as before, $ t^r $ times the right hand side of (\ref{righthandsideIPP}) is a linear combination of boundary terms of the form
\begin{eqnarray}
t^m \mbox{Im}\big( (\overline{P} - i \delta)^{-k} \big) \chi^{(j)} (0), \qquad 0 \leq m + k \leq r, \label{termedebord}
\end{eqnarray}
and of integrals of the form
\begin{eqnarray}
 \int_{0}^{\infty}  e^{\pm it \mu} \mu^l  \chi^{(j)} (\pm \mu^2) \mbox{Im} \left( ( \overline{P} \pm \mu^2 - i \delta )^{-k-1} \right)  d \mu , \label{integralereste}
\end{eqnarray}
with everywhere $$ 0 \leq j , k \leq r , \ \ l \geq 0, \qquad \mbox{ and } \ l\geq 1 \ \ \mbox{when } \ k =r . $$ 
By Theorem \ref{estimationprincipale}, the integrals (\ref{integralereste}) are uniformly bounded with respect to $ \delta $ and $t$ since the resolvents are bounded, except perhaps when $ k = r $ in which case they are at most of order $ |\mu|^{-1} $, but the latter is controlled by the term $ \mu^l $ with $ l \geq 1 $. To complete the proof, it suffices to show that the boundary terms
(\ref{termedebord}) are bounded with respect to $ \delta $ and $t$. This is clear if $ m = 0 $ since $ k \leq r $ then, and the resolvent to this power is bounded near the origin. It remains to show that (\ref{termedebord}) goes to zero  as $ \delta \rightarrow 0 $ if $ m \geq 1 $.   Indeed, by writing 
$$ \mbox{Im}\big( (\overline{P} - i \delta)^{-k} \big) = \frac{\delta}{2} \int_{-\pi/2}^{ \pi /2} (\overline{P} - \delta e^{i \theta})^{-k-1} e^{i \theta} d \theta , $$
and using that $ k \leq r - 1 $, we see that the limit is zero as $ \delta \rightarrow 0 $ since  we have a uniform bound for the resolvent inside the integral, since $ k+1 \leq r $. The result follows. \finpreuve

\bigskip







\noindent {\bf Proof of Theorem \ref{theoremeexponentiel}.} We study first the Schr\"odinger equation. We consider the second half of the partition of unit (\ref{hautefrequences}). 
Using Proposition \ref{deriveeexponentielle} and the same integration by part trick in the Stone formula as in the proof of Proposition \ref{etape1prop} (which is now simpler since we have no boundary term and no singularity), we see that if $ N $ is large enough, then
$$ \big| \big|\langle x \rangle^{-N} e^{-it\overline{P}} \varphi (h^2\overline{P}) \langle x \rangle^{-N} \big| \big|_{L^2 \rightarrow L^2} \lesssim \scal{t}^{-1} e^{C/h}, \qquad t \in \Ra, \ h \in (0,1] . $$
By interpolation between this bound and the trivial bound $ \big| \big| e^{-it\overline{P}} \varphi (h^2\overline{P}) \big| \big|_{L^2 \rightarrow L^2} \leq C $, we see that, for any $ \theta \in (0,1) $,
$$ \big| \big|\langle x \rangle^{-\theta N} e^{-it\overline{P}} \varphi (h^2\overline{P}) \langle x \rangle^{-\theta N} \big| \big|_{L^2 \rightarrow L^2} \lesssim \scal{t}^{-\theta} e^{C \theta /h}, \qquad t \in \Ra, \ h \in (0,1] . $$
Fix $ \nu > 0 $ and choose $ \theta  $ such that $ \nu = N \theta $, we then have the following alternative:
\begin{enumerate}
\item{in the region where $ e^{C/h} \scal{t}^{-\theta/2} \leq 1 $, we have
$$ \big| \big|\langle x \rangle^{- \nu} e^{-it\overline{P}} \varphi (h^2\overline{P}) \langle x \rangle^{- \nu} \big| \big|_{L^2 \rightarrow L^2} \lesssim \scal{t}^{-\theta/2} ,$$}
\item{in the region $ e^{C/h} \scal{t}^{-\theta/2} > 1 $, we have $ \log \scal{t} < 2 C / \theta h  $ so we obtain
$$ h^s \lesssim (1 + \log \scal{t})^{-s} , $$
and  have anyway the trivial bound
$$ \big| \big|\langle x \rangle^{-\nu} e^{-it\overline{P}} \varphi (h^2\overline{P}) \langle x \rangle^{-\nu} \big| \big|_{L^2 \rightarrow L^2} \leq C . $$ }
\end{enumerate}
This discussion shows that 
$$ h^s \big| \big|\langle x \rangle^{-\nu} e^{-it\overline{P}} \varphi (h^2\overline{P}) \langle x \rangle^{-\nu} \big| \big|_{L^2 \rightarrow L^2} \leq C (1 + \log \scal{t})^{-s}, \qquad t \in \Ra, \ \ h \in (0,1] , $$
and we conclude using Corollary \ref{corollairefinal} and Proposition \ref{etape1prop} to handle the low frequency part. More precisely, for the latter,  we interpolate between (\ref{borneainterpoler}) for $ \chi = \Phi_0 $ (see (\ref{hautefrequences})) and the trivial bound $ || \Phi_0 (\overline{P}) e^{-it\overline{P}}  ||_{L^2 \rightarrow L^2}  \leq C$  to be able to use the weight $ \scal{x}^{-\nu} $, in which case we still have a polynomial time decay rate hence a logarithmic one. 

The proof is completely similar for the wave and Klein-Gordon equations, using only the additional fact that, if $ \widetilde{\Phi} \Phi = \Phi $ and $ \widetilde{\Phi} \equiv 0 $ near $ 0 $, 
$$ \scal{x}^{-\nu} \widetilde{\Phi} (\overline{P}) (\overline{P}+m^2)^{-1/2} \scal{x}^{\nu}  $$
is a bounded operator from $ H^s $ to $ H^{s+1} $ for any fixed $ \nu \geq 0 $, $ m \geq 0 $ and $ s \in \Ra $.

\finpreuve

\bigskip

\noindent {\bf Proof of Theorem \ref{theoremenontrapping}.}  By the non trapping assumption, we have the semiclassical estimates (see \cite{Wang,RoENS})
$$ \big| \big| \scal{x}^{-\nu} \varphi (h^2\overline{P}) e^{-ith\overline{P}} \scal{x}^{-\nu} \big| \big|_{L^2 \rightarrow L^2} \leq C \scal{t}^{-s}, \qquad t \in \Ra , \ h \in (0,1], $$
provided that
$$  0 \leq s < \nu . $$
In the non semiclassical time scaling, this gives
$$ \big| \big| \scal{x}^{-\nu} \varphi (h^2\overline{P}) e^{-it\overline{P}} \scal{x}^{-\nu} \big| \big|_{L^2 \rightarrow L^2} \leq C \scal{t/h}^{-s} \leq C h^s \scal{t}^{-s} , $$
which, using the notation (\ref{notationnorme}), shows that
$$ h^{-s} e_{\nu}(h,t) \leq C \scal{t}^{-s}, \qquad t \in \Ra, \ h \in (0,1] . $$
Using Corollary \ref{corollairefinal}, we obtain
$$ \big| \big| \scal{x}^{-s} \Phi (\overline{P}) e^{-it\overline{P}} \scal{x}^{-s} \big| \big|_{H^{-s} \rightarrow L^2} \leq C \scal{t}^{-s}. $$
Since we may assume that $ 2(\bar{r}(d) +1)  < s < \nu $ the latter decays faster than $ \scal{t}^{-\bar{r}(d)} $ so the conclusion follows from Proposition \ref{etape1prop} and the remark thereafter. \finpreuve

\appendix

\section{Change of coordinates} \label{coordonneesadaptees}
\setcounter{equation}{0}

In this appendix, we recall how to choose a smooth diffeomorphism $ \chi : \Ra^d \rightarrow \Ra^d $ such that $ \chi^* G $ has determinant $ 1 $ outside a compact set and is still a long range perturbation of the euclidean metric. Recall that, if $ G = (G^{jk}) $ then,
\begin{eqnarray}
   \chi^* G = \big(\widetilde{G}^{jk}(y) \big) =  \ \mbox{Jac}_x(\chi)^{-1} (G^{jk}(x)) \mbox{Jac}_x(\chi)^{-1} , \qquad y = \chi (x) , \label{changementdecoordonneesmetrique}
\end{eqnarray}
where $ \mbox{Jac}_x (\chi) $ is the Jacobian matrix of $ \chi $ at $x$.
We shall show the following.
\begin{prop} Assume that $ 0 < \rho < 1 $. One can choose a smooth function $ \phi : \Ra^d \rightarrow \Ra $ such that, for some $ C > 0 $
\begin{enumerate}
\item{ for all $x \in \Ra^d $,
$$ C^{-1} \leq \phi (x) \leq C , $$}
\item{$ \phi -1 $ is a symbol of order $ - \rho $, ie 
$$ | \partial^{\alpha}_x ( \phi (x) - 1 ) | \leq C_{\alpha} \scal{x}^{- \rho-|\alpha|} , $$}
\item{The map $ \chi $ defined below is diffeomorphism from $ \Ra^d $ to $ \Ra^d $,
$$ \chi (x) = \phi(x) x , $$}
\item{for all $ |x| \geq C $,
\begin{eqnarray}
 \emph{det} (G^{jk}(x))^{1/2} & = & \phi (x)^{n-1} \left( \phi (x) +  x \cdot \nabla \phi (x) \right) , \label{transportcache} \\
& = &  \emph{det} \left( \emph{Jac}_x (\chi) \right) . \label{calculduJacobien}
\end{eqnarray}}
\item{The metric $ \widetilde{G} $ defined by (\ref{changementdecoordonneesmetrique}) is a long range perturbation of the euclidean metric, ie
$$ \big| \partial_y^{\alpha} \big( \widetilde{G}^{jk}(y) - \delta_{jk} \big) \big| \leq C_{\alpha} \scal{y}^{- \rho - | \alpha|} . $$}
\end{enumerate}
\end{prop}

\bigskip

\noindent {\it Proof.} We first solve (\ref{transportcache}) for $ |x| \geq R  $, for some $ R > 0 $ to be chosen. Since $ G $ is a long range perturbation of the euclidean metric, we have $ \mbox{det}(G)^{1/2} = 1 + \delta $, with $ \delta \in S^{-\rho} $.  By the change of unknown function $ \phi^n = \varphi $, (\ref{transportcache}) reads 
\begin{eqnarray}
 \varphi + \frac{x}{n} \cdot \nabla \varphi (x) = 1 + \delta (x), \qquad |x| \geq R . \label{transportreduit} 
\end{eqnarray}
This is a transport equation which is easily solved using polar coordinates; the solution which is equal to $1 $ on $ |x| = R $ is given by
\begin{eqnarray}
 \varphi (x) = 1 + n \int_1^{|x|/R} \delta (x/\tau) \frac{d\tau}{\tau^{n+1}}, \qquad |x| \geq R . \label{expressiontransport}
\end{eqnarray} 
The latter is well defined for any $ R > 0 $ and is clearly smooth. We shall choose $ R $ large enough to guarantee that $ \varphi $ is close enough to $ 1 $ and thus that $ \varphi^{1/n} $ is still smooth. Indeed, we have $ |x|/\tau \geq R $ on the interval of integration, thus
$$ |\varphi (x) - 1| \leq \sup_{|z|\geq R} |\delta(z)| , \qquad |x| \geq R ,$$
where the right hand side goes to zero as $ R \rightarrow \infty $. By multiplying $ \varphi $ by a smooth cutoff with values in $ [
0,1] $ which equals $1$ near infinity and $0$ near $ \{|x|\leq R \} $, we obtain a new function $ \varphi $ defined on $ \Ra^d $, such that $ |\varphi-1|\leq 1/2 $ everywhere which satisfies (\ref{transportreduit}) for a larger $ R $. Furthermore, by choosing $ R $ large enough, we may even assume that
\begin{eqnarray}
 \frac{1}{2} \leq \varphi (x) + \frac{x}{n} \cdot \nabla \varphi (x) \leq \frac{3}{2} . \label{constructiondiffeo}
\end{eqnarray}
We next check that $ \varphi - 1 \in S^{-\rho} $. Since the latter is a conditio at infinity, it is sufficient to consider the expression (\ref{expressiontransport}).
Using that $ \delta (z) \leq C |z|^{-\rho} $ we have,
$$ |\varphi (x) - 1| \leq  n \int_1^{\infty} |\delta (x/\tau)| \frac{d\tau}{\tau^{n+1}} \leq \frac{C}{|x|^{\rho}} \int_1^{\infty} \frac{\tau^{\rho}}{\tau^{n+1}} d \tau , $$
where the last integral is finite since $ \rho < 1 $. It is then not hard to check  that
$$ |\partial^{\alpha}_x (\varphi (x)-1)| \leq C_{\alpha} |x|^{-\rho-|\alpha|}, \qquad |x|> R , $$
by showing by induction that $ \partial^{\alpha}_x (\varphi (x)-1) $ is a linear combination of a symbol of order $ - |\alpha| - \rho $ and of terms of the form $ s_{\gamma} (x) \int_1^{|x|/R} (\partial^{\gamma} \delta )(x/\tau)
\tau^{-n-1-|\gamma|} d \tau$ with $s_{\gamma} $ a symbol of order $ | \gamma | - |\alpha|$, for $ |\gamma| \leq |\alpha| $.

  Setting $ \phi = \varphi^{1/n} $, we get a function satisfying the items $1$ and $2$, as well as (\ref{transportcache}). 
We now prove item 3. It is not hard to check that  $ \chi $ is a diffeomorphism if and only if, given $ \omega \in {\mathbb S}^{d-1} $, the map $ r \mapsto r \phi (r \omega) $ is diffeomorphism from $ \Ra^+ $ onto itself. The derivative of this function is
$$ \phi (r \omega) + r \omega \cdot \nabla \phi (r\omega) = \varphi (x)^{\frac{1}{n}-1}  \left( \varphi (x) + \frac{x}{n} \cdot \nabla \varphi (x) \right)_{x = r \omega} , $$
so the result follows from (\ref{constructiondiffeo}) and item 1. To prove (\ref{calculduJacobien}), one simply observes that
$$ \mbox{Jac}_x (\chi) =  \phi (x) I_d + \nabla \phi (x)^T x , $$
where the first matrix in the right hand side is scalar and the second matrix has rank one. The only possible non zero eignevalue of the latter is given by its trace which is $ x \cdot \nabla \phi (x) $. The caclulation of determinant is thus easy  and shows that (\ref{calculduJacobien}) coincides with the right hand side of (\ref{transportcache}). Finally, since $ \mbox{Jac}_x (\chi) - I_d $ is a matrix with entries in $ S^{-\rho} $ then so is the right hand side of (\ref{changementdecoordonneesmetrique}) and item 5 follows by a simple induction on $ |\alpha|$ using that $ \scal{x}^{-\rho} \approx \scal{\chi(x)}^{-\rho} $ and 
$ |\partial^{\alpha} (\chi(x)-x)| \lesssim \scal{x}^{-|\alpha|} $. \finpreuve


\begin{thebibliography}{99}



\bibitem{BTJFA} {\sc J.-M. Bouclet, N. Tzvetkov}, {\it On global Strichartz estimates for non
    trapping metrics}, J. Funct. Analysis 254 (2008) 1661-1682.
    
\bibitem{Bouc1} {\sc J.-M. Bouclet}, {\it Low frequency estimates for long range perturbations in divergence form}, CJM to appear.   

\bibitem{BoHa1} {\sc J.-F. Bony, D. H\"afner}, {\it The semilinear wave equation on asymptotically Euclidean manifolds}, CPDE to appear.

\bibitem{BoHa2} {\sc \name}, {\it Low frequency resolvent estimates for long range perturbations of the Euclidean Laplacian}, arXiv:0903.5531. 

\bibitem{Burq1} {\sc N. Burq},  {\it D\'ecroissance de l'\'energie locale de l'\'equation des ondes pour le probl\`eme ext\'erieur et absence de r\'esonance au voisinage du r\'eel}, Acta Math. 180 (1998) 1-29.  
    
\bibitem{Burq2} {\sc \name}, {\it Lower bounds for shape resonances widths of long range Schr\"odinger operators}, Amer. J. Math. 124, no. 4 (2002) 677-735.

\bibitem{stCardosoVodev} {\sc F. Cardoso, G. Vodev}, {\it Uniform estimates of the resolvent of the Laplace-Beltrami operator on infinite volume manifolds, II}, Ann. Henri Poincar\'e, 3, 673-691 (2002).

\bibitem{CardosoVodev2} {\sc \name}, {\it High frequency resolvent estimates and energy decay of solutions to the wave equation}, Canad. Math. Bull. vol 47 (4), (2004) 504-514.

\bibitem{Christianson} {\sc H. Christianson}, {\it Applications of Cutoff Resolvent Estimates to the Wave Equation},
 Math. Res. Lett. Vol. 16 (2009), no. 4, 577-590. 
 
\bibitem{Davies} {\sc E. B. Davies}, {\it Heat Kernels and Spectral Theory}, Cambridge Tracts in Math. Vol 92, Camb. Univ. Press (1989).
 

\bibitem{DiSj} {\sc M. Dimassi, J. Sj\"ostrand}, {\it  Spectral asymptotics in the semi-classical limit}, London Mathematical Society Lecture Note Series, 268. Cambridge University Press, Cambridge (1999).


\bibitem{DSS} {\sc R. Donninger. W. Schlag, A. Soffer}, {\it On pointwise decay for linear waves on a Schwarzschild black hole background}, preprint. 


\bibitem{FoSk} {\sc S. Fournais, E. Skibsted}, {\it Zero energy asymptotics of the resolvent for a class of slowly decaying potentials},  Math. Z.  248  (2004),  no. 3, 593-633.

\bibitem{GeMa} {\sc C. G\'erard, A. Martinez}, {\it Principe d'absorption limite pour des op\'erateurs de Schr\"odinger \`a longue port\'ee}, C.R. Acad. Sci. Paris 306, 121-123 (1988).

\bibitem{GuHa}  {\sc C. Guillarmou, A. Hassell}, {\it The resolvent at low energy and Riesz transform for Schr\"odinger operators on asymptotically conic manifolds, Part I}, 
Math. Ann. (2008) no 4, 859-896.

\bibitem{IsKi} {\sc H. Isozaki, I. Kitada}, {\it Modified wave operators with time independent
modifiers},  J. Fac. Sci. University of Tokyo, Section I A 32 (1985)
77-104.

\bibitem{IKRem} {\sc \name}, {\it A remark on  the micro-local resolvent estimates for two body Schr\"odinger operators}, Publ. RIMS, Kyoto Univ. {\bf 21} (1985), 889-910.

\bibitem{AJensen} {\sc A. Jensen}, {\it Propagation estimates for Schr\"odinger-type operators}, Trans. Amer. Math. Soc. 291, no. 1 (1985) 129-144.

\bibitem{JeMoPe} {\sc A. Jensen, E. Mourre, P. Perry}, {\it  Multiple commutator estimates and resolvent smoothness in quantum scattering theory}, Ann. IHP (A) Physique th\'eorique, 41 no. 2 (1984) 207-225.

\bibitem{KochTataru} {\sc H. Koch, D. Tataru}, {\it Carleman estimates and absence of embedded eigenvalues}, Commun. Math. Phys. 267, no. 2, 419-449 (2006).

\bibitem{LaPh} {\sc P. D. Lax, C. S.  Morawetz, R. S. Philips}, {\it Exponential decay of solutions of the wave equation in the exterior of a star-shaped obstacle},  Comm. Pure Appl. Math.  16  (1963) 477-486.

\bibitem{MeTa} {\sc D. J. Metcalfe, D. Tataru}, {\it Decay estimates for variable coefficient wave equations in exterior domains},
 Progress in Nonlinear Differential Equations and Their Applications, Vol. 78 (2009)  201-217. 



\bibitem{Naka1} {\sc S. Nakamura}, {\it Semiclassical resolvent estimates for the barrier top energy},
CPDE 16 (1991), no. 4-5, 873-883. 

\bibitem{Naka2} {\sc \name}, {\it Low energy asymptotics for Schrödinger operators with slowly decreasing potentials}, Comm. Math. Phys. 161(1), 63-76 (1994).

\bibitem{NoZw} {\sc S. Nonnenmacher, M. Zworski}, {\it Quantum decay rates in chaotic scattering}, Acta Math. Vol. 203, no. 2 (2009) 149-233. 

\bibitem{SSS} {\sc W. Schlag, A. Soffer and W. Staubach},   {\it Decay for the wave and Schr\"odinger evolutions on manifolds with conical ends. Part I},  
 Trans. Amer. Math. Soc.  362  (2010)  no. 1, 19-52.

\bibitem{ReedSimon} {\sc M. Reed, B. Simon}, {\it Methods of Modern Mathematical Physics I}, Academic Press (1980).

\bibitem{RoENS} {\sc D. Robert}, {\it Asymptotique de la phase de diffusion \`a haute \'energie pour des perturbations du second ordre du Laplacien}, 
Ann. Sci. de l'ENS, vol. 25, n. 2 (1992) 107-134.

\bibitem{RoTa} {\sc D. Robert, H. Tamura}, {\it Semiclassical estimates for resolvents and asymptotics for total scattering cross-sections}, Ann. Inst. H. Poincar\'e (phys. th\'eor.) 47, 415-442 (1987).

\bibitem{RodTao} {\sc I. Rodnianski, T. Tao}, {\it Longtime decay estimates for the Schr\"odinger equation on manifolds},  Mathematical aspects of nonlinear dispersive equations,  223-253, Ann. of Math. Stud., 163, Princeton Univ. Press, Princeton, NJ (2007). 


\bibitem{Tataru} {\sc D. Tataru}, {\it Parametrices and dispersive estimates for Schr\"odinger operators with variable coefficients}, Amer. J. Math.  130  (2008)  no. 3, 571-634.

\bibitem{Tataru2} {\sc \name}, {\it Local decay of waves on asymptotically flat stationary space-times}, preprint.

\bibitem{Vain} {\sc B. Vainberg}, {\it Asymptotic methods in equations of mathematical physics},  Gordon \etan Breach Science Publishers, New York, (1989).

\bibitem{VaWu} {\sc A. Vasy, J. Wunsch}, {\it Positive commutators at the bottom of the spectrum}, preprint.



\bibitem{VaZw} {\sc A. Vasy, M. Zworski}, {\it Semiclassical estimates in asymptotically Euclidean scattering},  
Comm. Math. Phys. 212 (2000), no. 1, 205-217.

\bibitem{Vodev} {\sc G. Vodev}, {\it Local energy decay of solutions to the wave equation for nontrapping metrics},
Ark. Mat. 42 (2004), 379-397.

\bibitem{Wang} {\sc X.P. Wang}, {\it Time-decay of scattering solutions and classical trajectories}, Ann. Inst. H. Poincar\'e (phys. th\'eor.) 47, 25-37 (1987).

\bibitem{Wang2} {\sc \name}, {\it Asymptotic expansion in time of the Schr\"odinger group on conical manifolds}, Ann. Inst. Fourier 56, 6, 1903-1945 (2006).

\bibitem{Yafa} {\sc D. Yafaev}, {\it The low energy scattering for slowly decreasing potentials}, Comm. Math. Phys. 85(2), 177--196 (1982).

\end{thebibliography}
\end{document}